\documentclass[english,11pt,a4paper]{article}

\usepackage{latexsym}         
\usepackage{amsmath}          
\usepackage{amssymb}          
\usepackage{amsxtra}          
\usepackage[latin1]{inputenc} 
\usepackage{eucal}            
\usepackage{longtable}        

\usepackage{exscale}          
\usepackage{bbm}              
\usepackage{theorem}          
\usepackage{paralist}         
\usepackage[arrow, matrix, curve]{xy}
\usepackage{stmaryrd}         
\usepackage{tikz}
\usepackage{tikz-cd}          
\usetikzlibrary{matrix, arrows}
\usepackage{pstricks}
\usepackage{enumerate}        
\usepackage{url}

\usepackage{leftidx}

\usepackage{a4wide}

\newcommand{\Doppel}[1]               {\mathcal{D}(#1)}
\newcommand{\opp}[1]             {{#1}^{\mathrm{op}}}
\newcommand{\ga}                 {\sideset{^{\scriptstyle{A}}}{}{\operatorname{g}}}
\newcommand{\xa}             {\sideset{^{\scriptstyle{A}}}{}{\operatorname{x}}}

\numberwithin{equation}{section}

\newcommand{\ques} [1] {\marginpar{\bf qu}{\bf {#1} ?}}
\newcommand{\chn}[1]{\marginpar{\bf changed}{\bf {#1}}}
\newcommand{\note} [1] {\marginpar{\bf note}{\bf {#1}}}
\newcommand{\todo} [1] {\marginpar{\bf to do}{\bf {#1}}}
\newcommand{\prob} [1] {\marginpar{\bf problem}{\bf {#1}}}

\newcommand{\I}          {\mathrm{i}}

\newcommand{\C}        {\mathbb{C}}

\newcommand{\N}        {\mathbb{N}}

\newcommand{\Z}        {\mathbb{Z}}

\newcommand{\Q}        {\mathbb{Q}}

\newcommand{\R}        {\mathbb{R}}

\newcommand{\E}          {\mathrm{e}}

\newcommand{\D}          {\operatorname{\mathrm{d}}}

\newcommand{\cc}[1]      {\overline{{#1}}}

\newcommand{\sign}       {\operatorname{\mathrm{sign}}}

\newcommand{\RE}         {\mathsf{Re}}

\newcommand{\IM}         {\mathsf{Im}}

\newcommand{\Unit}       {\mathrm{Id}}

\newcommand{\cl}         {\mathrm{cl}}

\newcommand{\const}      {\mathit{const}}

\newcommand{\at}[1]      {\big|_{#1}}

\newcommand{\At}[1]      {\Big|_{#1}}

\newcommand{\argument}   {\,\cdot\,}

\newcommand{\pr}             {\mathrm{pr}}

\newcommand{\El}[5] {\sideset{^{\scriptscriptstyle{#1}}_{\scriptscriptstyle{#2}}}{^{\scriptscriptstyle{#4}}_{\scriptscriptstyle{#5}}}{\operatorname{#3}}}

\newcommand{\tensor}[1][{}]   {\mathbin{\otimes_{\scriptscriptstyle{#1}}}}

\newcommand{\boxtensor}[1][{}] {\mathbin{\boxtimes_{\scriptscriptstyle{#1}}}}

\newcommand{\alg}[1]            {\mathcal{#1}}
\newcommand{\ring}[1]            {\mathcal{#1}}

\newcommand{\swe}[1]    {{\scriptscriptstyle{(#1)}}}

\newcommand{\we}[1]    {{\scriptscriptstyle{#1}}}

\newcommand{\id}         {\operatorname{\mathsf{id}}}

\newcommand{\ad}         {\operatorname{\mathrm{ad}}}

\newcommand{\Ad}         {\operatorname{\mathrm{Ad}}}

\newcommand{\Mat}         {\operatorname{\mathrm{Mat}}}

\newcommand{\Conj}       {\operatorname{\mathrm{Conj}}}
\newcommand{\tra}[1]       {{#1}^{\operatorname{\mathrm{t}}}}
\newcommand{\rr}       {\rightarrow}
\newcommand{\mto}  {\mapsto}
\newcommand{\lrr}          {\longrightarrow}
\newcommand{\lmto}      {\longmapsto}
\newcommand{\act}  {\triangleright}

\newcommand{\oppact} {\opp{\triangleright}}
\newcommand{\ractopp} {\opp{\triangleleft}}

\newcommand{\ract} {\triangleleft}

\newcommand{\Cat}[1]         {\operatorname{\mathcal{#1}}}
\newcommand{\op}[1]           {\Cat{#1}^{\operatorname{\mathrm{op}}}}
\newcommand{\rev}[1]           {\Cat{#1}^{\operatorname{\mathrm{rev}}}}

\newcommand{\Hom}        {\operatorname{\mathsf{Hom}}}

\newcommand{\Homo}        {\operatorname{\mathsf{H}}}

\newcommand{\Chain}        {\operatorname{\mathsf{C}}}

\newcommand{\End}        {\operatorname{\mathsf{End}}}
\newcommand{\Gl}        {\operatorname{\mathsf{Gl}}}
\newcommand{\Aut}        {\operatorname{\mathsf{Aut}}}
\newcommand{\Iso}        {\operatorname{\mathsf{Iso}}}
\newcommand{\Der}        {\operatorname{\mathsf{Der}}}
\newcommand{\Obj}        {\operatorname{\mathsf{Obj}}}
\newcommand{\Mor}        {\operatorname{\mathsf{Mor}}}
\newcommand{\InnAut}     {\operatorname{\mathsf{InnAut}}}
\newcommand{\OutAut}     {\operatorname{\mathsf{OutAut}}}
\newcommand{\InnDer}     {\operatorname{\mathsf{InnDer}}}
\newcommand{\OutDer}     {\operatorname{\mathsf{OutDer}}}

\newcommand{\K}     {\operatorname{\mathsf{K}}}

\newcommand{\tr}         {\operatorname{\mathsf{tr}}}

\newcommand{\Vect}         {\operatorname{\mathsf{Vect}}}

\newcommand{\zentrum}    {\mathcal{Z}}

\newcommand{\Mod}      {\operatorname{\mathsf{Mod}}_{\Cat{C}}}
\newcommand{\ModC}      {\operatorname{\mathsf{Mod}}{(\Cat{C})}}
\newcommand{\ModCt}      {{\mathsf{Mod}}^{\Theta}{(\Cat{C})}}

\newcommand{\RMod}      {\operatorname{\mathsf{R-Mod}}}

\newcommand{\LMod}      {\operatorname{\mathsf{L-Mod}}}

\newcommand{\Rep}      {\operatorname{\mathsf{Rep}}}

\newcommand{\Cats}      {\operatorname{\mathsf{Cat}}}

\newcommand{\Comod} {\operatorname{\mathsf{Comod}}}

\newcommand{\ev}[1]   {\operatorname{\mathsf{ev}}_{#1}}
\newcommand{\coev}[1]   {\operatorname{\mathsf{coev}}_{#1}}
\newcommand{\evp}[1]   {\operatorname{\mathsf{ev}}_{#1}^{\prime}}
\newcommand{\coevp}[1]   {\operatorname{\mathsf{coev}}_{#1}^{\prime}}

\newcommand{\FPdim}      {\operatorname{\mathsf{FPdim}}}

\newcommand{\Bimod}[5]{\sideset{^{\scriptscriptstyle{#1}}_{\scriptscriptstyle{#2}}}{^{\scriptscriptstyle{#4}}_{\scriptscriptstyle{#5}}}{\operatorname{#3}}}

\newcommand{\DF}  {\Bimod{}{}{\Cat{D}}{F}{}}
\newcommand{\DG}  {\Bimod{}{}{\Cat{D}}{G}{}}
\newcommand{\CG}  {\Bimod{}{}{\Cat{C}}{G}{}}
\newcommand{\DFG}  {\Bimod{}{}{\Cat{D}}{FG}{}}
\newcommand{\bimCD} {(\Cat{C}, \Cat{D})}
\newcommand{\bimDC} {(\Cat{D}, \Cat{C})}

\newcommand{\CMstar}  {\Bimod{}{}{\Cat{C}}{*}{\Cat{M}}}
\newcommand{\CMphi}  {\Bimod{\varphi}{\Cat{C}}{\Cat{M}}{}{}}
\newcommand{\CNphi}  {\Bimod{\varphi}{\Cat{C}}{\Cat{N}}{}{}}
\newcommand{\CEphi}  {\Bimod{\varphi}{\Cat{C}}{\Cat{E}}{}{}}
\newcommand{\CCphi}  {\Bimod{\varphi}{\Cat{C}}{\Cat{C}}{}{}}
\newcommand{\MCphi}  {\Bimod{}{}{\Cat{M}}{\Cat{C}}{\varphi}}

\newcommand{\CM}  {\Bimod{}{\Cat{C}}{\Cat{M}}{}{}}
\newcommand{\CMi}  {\Bimod{}{\Cat{C}}{\Cat{M}}{i}{}}
\newcommand{\CMone}  {\Bimod{}{\Cat{C}}{\Cat{M}}{1}{}}
\newcommand{\CMn}  {\Bimod{}{\Cat{C}}{\Cat{M}}{n}{}}
\newcommand{\MC}  {\Bimod{}{}{\Cat{M}}{}{\Cat{C}}}
\newcommand{\CN}  {\Bimod{}{\Cat{C}}{\Cat{N}}{}{}}
\newcommand{\NC}  {\Bimod{}{}{\Cat{N}}{}{\Cat{C}}}
\newcommand{\CE}  {\Bimod{}{\Cat{C}}{\Cat{E}}{}{}}

\newcommand{\DM}  {\Bimod{}{\Cat{D}}{\Cat{M}}{}{}}
\newcommand{\MD}  {\Bimod{}{}{\Cat{M}}{}{\Cat{D}}}
\newcommand{\DN}  {\Bimod{}{\Cat{D}}{\Cat{N}}{}{}}
\newcommand{\ND}  {\Bimod{}{}{\Cat{N}}{}{\Cat{D}}}

\newcommand{\CMD}  {\Bimod{}{\Cat{C}}{\Cat{M}}{}{\Cat{D}}}
\newcommand{\CND}  {\Bimod{}{\Cat{C}}{\Cat{N}}{}{\Cat{D}}}
\newcommand{\DMC}  {\Bimod{}{\Cat{D}}{\Cat{M}}{}{\Cat{C}}}
\newcommand{\DNC}  {\Bimod{}{\Cat{D}}{\Cat{N}}{}{\Cat{C}}}
\newcommand{\CNF}  {\Bimod{}{\Cat{C}}{\Cat{N}}{}{\Cat{F}}}
\newcommand{\CNE}  {\Bimod{}{\Cat{C}}{\Cat{N}}{}{\Cat{E}}}

\newcommand{\DNF}  {\Bimod{}{\Cat{C}}{\Cat{N}}{}{\Cat{F}}}

\newcommand{\CYF}  {\Bimod{}{\Cat{C}}{\Cat{Y}}{}{\Cat{F}}}

\newcommand{\DDD}  {\Bimod{}{\Cat{D}}{\Cat{D}}{}{\Cat{D}}}
\newcommand{\CCC}  {\Bimod{}{\Cat{C}}{\Cat{C}}{}{\Cat{C}}}
\newcommand{\CC}  {\Bimod{}{\Cat{C}}{\Cat{C}}{}{}}

\newcommand{\CMDop}  {\Bimod{}{\Cat{C}}{\Cat{M}}{\mathrm{op}}{\Cat{D}}}
\newcommand{\DMCop}  {\Bimod{}{\Cat{D}}{\Cat{M}}{\mathrm{op}}{\Cat{C}}}

\newcommand{\MCop}  {\Bimod{}{\Cat{}}{\Cat{M}}{\mathrm{op}}{\Cat{C}}}
\newcommand{\CMop}  {\Bimod{}{\Cat{C}}{\Cat{M}}{\mathrm{op}}{\Cat{}}}

\newcommand{\NCop}  {\Bimod{}{\Cat{}}{\Cat{N}}{\mathrm{op}}{\Cat{C}}}
\newcommand{\CNop}  {\Bimod{}{\Cat{C}}{\Cat{N}}{\mathrm{op}}{\Cat{}}}

\newcommand{\Lmult}[2]   {\Bimod{#1}{}{L}{}{#2}}
\newcommand{\Lmultrr}[2]   {\Bimod{#1}{}{L}{rr}{#2}}

\newcommand{\lmult}[2]   {\Bimod{#1}{}{l}{}{#2}}

\newcommand{\IP}[4]{{\,}_{\scriptscriptstyle{#2}\!\!}\left\langle{{#1}}\right\rangle^{\scriptscriptstyle{#3}}_{\scriptscriptstyle{#4}}}

\newcommand{\IPs}[4]{{\,}_{\scriptscriptstyle{#2}\!\!}\left\langle{{#1}}\right\rangle^{{#3}}_{\scriptscriptstyle{#4}}}
\newcommand{\IPl}[4]{{\,}_{\scriptscriptstyle{#2}\!\!}^{*}\left\langle{{#1}}\right\rangle^{\scriptscriptstyle{#3}}_{\scriptscriptstyle{#4}}}

\newcommand{\icmM}[1]    {\underline{\mathsf{Hom}}({#1})^{\scriptscriptstyle{\mathcal{M}}}}

\newcommand{\icm}[1]      {\underline{\mathsf{Hom}}({#1})}

\newcommand{\imd}[1]    {\IP{{#1}}{}{\Cat{M}}{\Cat{D}}}

\newcommand{\icmr}[1]    {\IPs{{#1}}{\mathcal{C}}{r}{}}
\newcommand{\ic}[1]    {\IP{{#1}}{}{}{}}

\newcommand{\istarcm}[1]   {\IPl{{#1}}{\mathcal{C}}{}{}}

\newcommand{\icmod}[1]    {\IP{{#1}}{\mathcal{C}}{\mathsf{Mod}(A)}{}}

\newcommand{\icmstar}[1]    {\IPs{{#1}}{\mathcal{C}}{*}{}}
\newcommand{\icstar}[1]    {\IPs{{#1}}{}{*}{}}

\newcommand{\idm}[1]    {\IP{{#1}}{\Cat{D}}{\Cat{M}}{}}
\newcommand{\imc}[1]    {\IP{{#1}}{}{\mathcal{M}}{\Cat{C}}}

\newcommand{\icn}[1]    {\IP{{#1}}{\Cat{C}}{\Cat{N}}{}}
\newcommand{\ine}[1]    {\IP{{#1}}{}{\mathcal{N}}{\Cat{E}}}

\newcommand{\idmn}[1]    {\IP{{#1}}{\Cat{D}}{\Cat{M} \Box \Cat{N}}{}}
\newcommand{\imne}[1]    {\IP{{#1}}{}{\Cat{M} \Box \Cat{N}}{\Cat{E}}}
\newcommand{\iAm}[1]    {\IP{{#1}}{}{A}{}}
\newcommand{\Funl}[2]  {\operatorname{\mathrm{Fun}}_{\scriptscriptstyle{#1}}{(#2)}}
\newcommand{\FunDF}  {\Funl{\Cat{D}}{\DF,\DF}}

\newcommand{\Funr}[2]{\operatorname{\mathrm{Fun}}{(#2)}_{\scriptscriptstyle{#1}}}

\newcommand{\UFunl}[2]  {\operatorname{\underline{\mathrm{Fun}}}_{\scriptscriptstyle{#1}}{(#2)}}
\newcommand{\UFunDF}  {\UFunl{\Cat{D}}{\DF,\DF}}

\newcommand{\UFunr}[2]  {\operatorname{\underline{\mathrm{Fun}}}{(#2)}_{\scriptscriptstyle{#1}}}

\newcommand{\fmu}[1] {f^{\mu}_{#1}}
\newcommand{\fnu}[1] {f^{\nu}_{#1}}
\newcommand{\frho}[1] {f^{\rho}_{#1}}

\newcommand{\tenrev} {\operatorname{\otimes}^{\scriptscriptstyle{\mathrm{rev}}}}
\newcommand{\btD}    {\boxtensor[\Cat{D}]}
\newcommand{\btC}    {\boxtensor[\Cat{C}]}

\newcommand{\unD}  {1_{\Cat{D}}}
\newcommand{\unC}  {1_{\Cat{C}}}

\newcommand{\BimCat}   {\underline{\underline{\mathsf{BimCat}}}}  
\newcommand{\Ps}          {\underline{\underline{\mathsf{Ps}}}}
\newcommand{\Gray}   {\bold{\mathcal{G}}}
\newcommand{\opped}[1]   {{#1}^{\mathrm{op}_{\scriptscriptstyle{12}}}}

\newcommand{\corres}{\mathrel{\widehat{=}}}

\theoremheaderfont{\normalfont\bfseries}
\theorembodyfont{\itshape}
\newtheorem{lemma}{Lemma}[section]
\newtheorem{proposition}[lemma]{Proposition}
\newtheorem{theorem}[lemma]{Theorem}
\newtheorem{corollary}[lemma]{Corollary}
\newtheorem{definition}[lemma]{Definition}

\theorembodyfont{\rmfamily}

\newtheorem{example}[lemma]{Example}
\newtheorem{remark}[lemma]{Remark}

\newenvironment{theorem-n}[1][Theorem]{\begin{trivlist}
  \item[\hskip \labelsep {\bfseries #1}]}{\end{trivlist}}

\newcommand{\fixspacingbeforelist}{\mbox{}}
\makeatletter
\newcommand{\rmnum}[1]{\romannumeral #1}
\newcommand{\Rmnum}[1]{\expandafter\@slowromancap\romannumeral #1@}
\makeatother

\renewcommand{\labelenumi}{{\itshape\roman{enumi})}}
\renewcommand{\labelitemi}{$-$}

\makeatletter
\renewcommand{\theenumi}{\roman{enumi}}
\makeatother
\newcommand{\aditem}[1]  {\noindent\textit{ad \ref{#1}):}}
\newcommand{\refitem}[1] {~\textit{\ref{#1})}}

\makeatletter
\newcommand\qedsymbol{\hbox{$\boxempty$}}
\newcommand\qed{\relax\ifmmode\boxempty\else
  {\unskip\nobreak\hfil\penalty50\hskip1em\null\nobreak\hfil\qedsymbol
    \parfillskip=\z@\finalhyphendemerits=0\endgraf}\fi}
\makeatother
\newenvironment{proof}[1][{}]{\par\noindent \emph{Proof}{#1}: }{\qed}

\newenvironment{theoremlist}{\begin{enumerate}}{\end{enumerate}}
\newenvironment{remarklist}{\begin{enumerate}}{\end{enumerate}}

\newenvironment{lemmalist}{\begin{enumerate}}{\end{enumerate}}
\newenvironment{propositionlist}{\begin{enumerate}}{\end{enumerate}}
\newenvironment{definitionlist}{\begin{enumerate}}{\end{enumerate}}

\title{Traces on Module Categories over Fusion Categories}
\author{
  \textbf{Gregor Schaumann}\thanks{email:
    gregor.schaumann@math.uni-erlangen.de}
  \\[0.1cm]
  Department Mathematik \\
  Friedrich-Alexander Universit\"at Erlangen-N\"urnberg\\
  Cauerstra{\ss}e 11\\
  91058 Erlangen \\
  Germany}

\begin{document}
\maketitle
\abstract{We  consider  traces on module categories over pivotal fusion categories which are compatible with the module structure. 
It is shown that such module traces characterise the  Morita classes of special haploid symmetric Frobenius algebras. 
Moreover, they are unique up to a scale factor and they equip the dual category with a pivotal structure. 
This implies that for each  pivotal structure on 
a fusion category over $\C$ there exists a conjugate pivotal structure defined by the canonical module trace.
}

\section{Introduction}

Fusion categories  exhibit a rich mathematical structure, see for example \cite{ENOfus, Mue1}. They 
have important applications in  3-dimensional topological field theory \cite{BakKir, Resh}, in particular in the study of invariants of 3-manifolds \cite{BarWes, Tur},  and  in
rational conformal field theory, see \cite{MooreSei}, \cite{SchwPart} and subsequent work.
The construction  in conformal field theory initiated in  \cite{SchwPart} requires as its starting point a special haploid symmetric Frobenius object in a modular fusion category, but it depends only on the 
Morita class of that algebra. It is known \cite{Ostrik} that Morita classes of algebras in fusion categories are described by equivalence classes of module categories.

In this article we provide a description of the Morita classes of special  haploid symmetric Frobenius algebras in pivotal fusion categories over $\C$ in terms of module categories with module traces.
 A module trace is a trace on a module category, i.e.~a collection of  symmetric and non-degenerate linear maps from the endomorphism spaces of  objects to $\C$, that is compatible with the module structure. 
As a main result we prove  the following:

\begin{theorem-n}
\emph{Let $\Cat{C}$ be a pivotal fusion category. The following structures are equivalent:}
\begin{enumerate} 
\item 
  \label{item:module-cat}
  \emph{An indecomposable module category $\CM$ with module trace.}
\item \label{item:C-bal}\emph{An indecomposable module category $\CM$ together with a $\Cat{C}$-balanced natural isomorphism between  $\Hom(n,m)$ and the dual space of  $\Hom(m,n)$, for each pair of objects $m,n\in  \Cat{M}$. 
\item \label{item:Frob-Morita} A Morita class of a special  haploid symmetric Frobenius algebra in $\Cat{C}$.} 
\end{enumerate}
\end{theorem-n}
The equivalence of \refitem{item:module-cat} and \refitem{item:C-bal} implies that module traces on indecomposable module 
categories are unique up to a constant factor and equip the dual fusion category with a pivotal structure.
When applied to the particular case of 
$\Cat{C}$ considered as a left module category over itself, we obtain the following result.
\begin{theorem-n}
\emph{For each pivotal structure $a$ on a fusion category $\Cat{C}$ over  $\C$ there exists a conjugate pivotal structure $\cc{a}$ 
such that the left dimensions of objects with respect to $\cc{a}$ are complex conjugate to the left dimensions 
with respect to $a$}. 
\end{theorem-n}
We show how this result is related to the existence of a natural monoidal isomorphism of the identity and the quadruple dual functor for fusion categories from \cite{ENOfus}.

We give an explicit description of module traces in terms of a matrix equation that 
    provides  a reduction of the problem of solving a quadratic equation for  algebras (the Frobenius property) to a linear equation for the module category.
This implies in particular  that the quantum dimensions of special haploid symmetric Frobenius algebras in pivotal fusion categories are positive real numbers and shows 
that all module categories over pseudo-unitary fusion categories admit a module trace.
We extend the graphical calculus for tensor categories to module categories and give a graphical description 
of the Frobenius algebra obtained from a module category with module trace. 

In \cite{DavWit} it is shown that indecomposable module categories over a fusion category $\Cat{C}$ are classified by Lagrangian algebras in the Drinfeld center $\mathcal{Z}(\Cat{C})$. It remains to 
interpret our results in terms of this classification. 

A possible application of our results is  to modify the construction in \cite{SchwPart} in such a way that it depends  only on a  module category with module traces over a modular fusion category and involves no further choices.
In such a construction it should be possible to incorporate module functors and module natural transformations as well and interpret them in physical terms, see  \cite[Sec. 3]{SchwCatCor}, \cite{KapSau}, \cite{FuSchwVal} for  a possible 
interpretation.

The paper is structured as follows.
In Section \ref{sec:preliminaries} we summarise the relevant background  about  fusion categories, algebra objects and module categories.
In Section \ref{sec:module-traces} we first develop a graphical notation for module categories which gives rise to  a diagrammatic  description of
 the algebra structure of inner hom objects. 
Next we introduce module traces and demonstrate in examples that the existence of a module trace for a given module category depends on the choice of pivotal structure for the fusion category.
In v3 of this article we clarify the existence of matched pairs of pivotal structure and module trace 
using ideas provided by Pavel Etingof. 
In Section \ref{sec:module-traces-dual} we give a description of module traces in terms of $\Cat{C}$-balanced natural isomorphisms  and prove that module traces on indecomposable module categories are unique up to scaling.
This description of module traces yields a module natural isomorphism between a  module functor and its double adjoint functor. In the application to a
pivotal fusion category as a module category over itself, this  leads to the existence of 
 conjugate pivotal structures for pivotal fusion categories. We provide a graphical derivation of a monoidal natural isomorphism of the identity functor to the quadruple dual functor for fusion categories and 
show that this yields an alternative definition of the conjugate pivotal structure. 
In Section \ref{sec:existence-problem-as} we demonstrate that the existence of a module trace can be reduced to a matrix equation and  discuss the example of pseudo-unitary fusion categories. 
As a consequence of these results we obtain a new criterion to decide whether a pivotal structure is spherical in terms of module categories.
In  Section \ref{subsec:Frobenius-Algebras}  we prove that module traces characterise equivalence classes of special haploid symmetric Frobenius algebras.

\section{Preliminaries}
\label{sec:preliminaries}
\subsection{Fusion Categories and Algebra Objects}
In this section we  summarise  the relevant background and fix our notation. 
All categories are assumed to be abelian and moreover locally finite over $\C$, i.e.~the isomorphism classes of objects form a set, all $\Hom_{}$-spaces are finite dimensional and every object has finite length.
All functors and natural transformations are assumed to be additive.
\begin{definition}{\rm \cite{EGNObook}}
  A tensor category $\Cat{C}$ is a monoidal category with rigidity and simple unit $1 \in \Cat{C}$ such that the monoidal structure is bilinear on morphisms. A finite tensor category is a tensor category with 
  finitely many simple objects up to isomorphism. A fusion category is a semisimple finite tensor category.
\end{definition}
 Without loss of generality we will work with strict monoidal categories (see e.g.~\cite{BakKir}).
Rigidity means that each object $c \in \Cat{C}$ has a right dual $c^{*}$ with duality morphisms 
$\ev{c}: c^{*} \otimes c \rr 1$,
$\coev{c}: 1 \rr c \otimes c^{*}$ and a left dual $\leftidx{^*}{c}{}$ with $\evp{c}: c \otimes \leftidx{^*}{c}{} \rr 1$ and $\coevp{c}: 1 \rr \leftidx{^*}{c}{} \otimes c $, 
such that the rigidity axioms are satisfied, see Appendix \ref{subsection:Graphical-calculus}, equation (\ref{rigidity}). Right and left 
duals are unique up to a unique isomorphism.
In a rigid tensor category there is a canonical 
natural isomorphism $c \simeq \leftidx{^*}{(c^{*})}{} \simeq ( \leftidx{^*}{c}{} )^{*}$ for each object $c \in \Cat{C}$ and  we will therefore identify these objects  in the sequel.

The functor $ (.)^{**}$ has a canonical structure of a tensor functor. 
A pivotal structure for $\Cat{C}$ is a monoidal natural  isomorphism $a:\id_{\Cat{C}}\rr  (.)^{**}  $.
 In particular, a pivotal structure allows one to define the left trace of a morphism $f \in \End(c)$ as 
\begin{equation}
  \label{eq:right-trace}
  \tr_{c}^{L}(f)= \ev{c}\circ (a_{ \leftidx{^*}{c}{} } \otimes  f) \circ \coevp{c} \in \End(1)\simeq \C
\end{equation}
and for each object $c$ the quantum dimension 
$\tr_{c}^{L}(\id_{c})=\dim^{\Cat{C}}(c)$. The right trace of a morphism is defined analogously and a pivotal structure is called spherical if the  left traces and right traces agree for all morphisms.
Throughout this paper $\Cat{C}$ denotes a pivotal fusion category unless stated otherwise.
We use the well-established graphical calculus for tensor categories,  see Appendix \ref{subsection:Graphical-calculus} for relevant definitions and conventions.

\paragraph{Algebra Objects}
\begin{definition}
  An algebra (object)  in a tensor category $\Cat{C}$ is an object $A \in \Cat{C}$ together with a multiplication morphism 
  $\mu: A\otimes A \rightarrow A$,   and a unit morphism $\eta :1 \rightarrow A$, represented by the diagrams
  \begin{equation}
    \mu \corres 
\ifx\du\undefined
  \newlength{\du}
\fi
\setlength{\du}{10\unitlength}

 \; ,
  \end{equation}
  such that the   associativity and unit constraints hold:
  \begin{equation}
    \label{eq:Algebra-definition}
\ifx\du\undefined
  \newlength{\du}
\fi
\setlength{\du}{10\unitlength}
%
\;.
  \end{equation}
  An algebra $A$ in $\Cat{C}$ is called haploid if $\Hom_{\Cat{C}}(1,A) \simeq \C$ as a vector space. 
\end{definition}
There is the obvious definition of morphisms of algebras. An algebra is called indecomposable if it is not isomorphic 
to a direct sum of two non-trivial algebras.
As we will always work with just one algebra at a time, we omit the labels on the lines representing 
the algebra object.
Given an algebra in $\Cat{C}$, we can consider modules over this algebra in $\Cat{C}$. 
\begin{definition}
  A right module over an algebra $A$ in a tensor category $\Cat{C}$ is an object $M \in \Cat{C}$ together with an action morphism 
  \begin{equation}
    \rho: M \otimes A \rightarrow M  \quad \corres \quad  
\ifx\du\undefined
  \newlength{\du}
\fi
\setlength{\du}{10\unitlength}

\;,
  \end{equation}
  such that the following  equations hold:
  \begin{equation}
\ifx\du\undefined
  \newlength{\du}
\fi
\setlength{\du}{10\unitlength}
%
\;.
  \end{equation}

  An intertwiner between two right modules $(M, \rho)$ and $(N, \chi)$ over $A$ is a morphism $\phi: M \rightarrow N$ in $\Cat{C}$ which satisfies
  \begin{equation}
\ifx\du\undefined
  \newlength{\du}
\fi
\setlength{\du}{10\unitlength}

\;.
  \end{equation}
 There are analogous definitions for left modules.
  The subspace of $\Hom_{\Cat{C}}(M,N)$
  consisting of  the intertwiners is denoted by  $\Hom_{A}(M,N)$. 
\end{definition}

It is clear (see e.g.~\cite{Ostrik}) that for an algebra $A$, a right module  $(M, \rho)$ over $A$ and an object $c \in \Cat{C}$, the object $c \otimes M$ is also a right module over $A$ with action morphism 
  \begin{equation}\label{eq:tensor-module}
    \id_c \otimes \rho :c \otimes M \otimes A \rightarrow c \otimes M,
  \end{equation}
and that each   morphism $\phi: c \rightarrow d$ in $\Cat{C}$ yields an intertwiner $\phi \otimes \id_M: c \otimes M \rightarrow d \otimes M$.

\begin{definition}{\rm \cite{SchwCat}}
  An coalgebra (object)  in a tensor category $\Cat{C}$ is an object $C \in \Cat{C}$ together with a comultiplication morphism 
  \begin{equation}
    \Delta:  C \rr C \otimes C  \quad \corres \quad  
\ifx\du\undefined
  \newlength{\du}
\fi
\setlength{\du}{10\unitlength}

\;,
  \end{equation}
  such that the coassociativity and counit constraints hold:
  \begin{equation}
\ifx\du\undefined
  \newlength{\du}
\fi
\setlength{\du}{10\unitlength}
%
\;.
  \end{equation} 
\end{definition}

\begin{definition}{\rm\cite{DavWit,SchwCat}}
  \label{definition:Frobenius-algebra}
 Let $\Cat{C}$ be a tensor category.
  \begin{definitionlist}
  \item A separable algebra $A \in \Cat{C}$ is an algebra  $(A, \mu, \eta)$  for which there exists a morphism  $\Delta: A \rr A \otimes A$ with $\mu \circ \Delta=\id_{A}$ and
    \begin{equation}
      \label{eq:splitting-multiplication}
      \Delta \circ \mu = (\mu \otimes \id_{A}) \circ (\id_{A} \otimes \Delta)= (\id_{A} \otimes \mu) \circ (\Delta \otimes \id_{A}). 
    \end{equation}
  \item  A Frobenius algebra in $\Cat{C}$ is an algebra $(A, \mu, \eta)$ that is also a coalgebra 
    with structures $\epsilon: A \rr 1$ and $\Delta: A\rr A \otimes A$, such that (\ref{eq:splitting-multiplication}) is satisfied.
  \end{definitionlist}
\end{definition}
In graphical notation relation (\ref{eq:splitting-multiplication}) reads :
\begin{equation}
  \label{eq:Frob-equ-graphical}
\ifx\du\undefined
  \newlength{\du}
\fi
\setlength{\du}{10\unitlength}

\;. 
\end{equation}
\begin{lemma}{\rm \cite[Prop. 2.7]{DavWit}}
  Consider an algebra  $(A, \mu, \eta)$ in a fusion category $\Cat{C}$. 
  Then the category $\Mod(A)$ is semisimple if and only if $A$ is separable.
\end{lemma}
Frobenius algebras with the following additional properties are particularly important in applications to conformal field theory 
\cite{SchwPart}. Here and in the following we adopt the convention, that an unlabelled small box in a diagram 
represents an isomorphism obtained from the pivotal structure.

\begin{definition}{\rm\cite{SchwCat}}
  A Frobenius algebra $A$ in $\Cat{C}$ is called
  \begin{definitionlist}
  \item  special if there exist $\beta_{1},\beta_{A} \in \C^{\times}$ such that
    \begin{equation}
      \label{special-diagram}
\ifx\du\undefined
  \newlength{\du}
\fi
\setlength{\du}{10\unitlength}
  \;.
    \end{equation}
  \end{definitionlist} 
\end{definition}
 Condition \refitem{item:symmetric-frob} makes sense for any algebra $A$ with a morphism $\epsilon \in \Hom_{\Cat{C}}(A,1)$.
\begin{lemma}{\rm\cite{FFRS}}
  \label{lemma:A-special-symmetric-dim-neq-zero}
  Let $A$ be a special symmetric Frobenius algebra in $\Cat{C}$. Then $\dim^{\Cat{C}}(A) =
  \beta_{1} \beta_{A} \neq 0$. We can normalise $\epsilon$ and $\Delta$ such that $\beta_{1}=\dim^{\Cat{C}}(A)$
  and $\beta_{A}=1$.
\end{lemma}

\begin{lemma}{\rm \cite{SchwPart}}
  If an algebra $A$ is haploid and has dimension $\dim^{\Cat{C}}(A) \neq 0$ \footnote{ In the proof 
  \cite[Cor. 3.10]{SchwPart} the assumption $\dim^{\Cat{C}}(A) \neq 0$ is implicitly present. We thank I. Runkel for this information.}
  \label{lemma:special-haploid-then-symmetric}, then it is symmetric for any choice of $\epsilon \in \Hom_{\Cat{C}}(A,1)$. 
\end{lemma}
Let  $\Cat{C}$ be a pivotal fusion category. 
The left dual $\leftidx{^*}{M}{}$ of a right $A$-module $(M, \rho)$ inherits a  canonical left $A$-module structure  
\begin{equation}
  \label{equation:dual-action}
  \rho_{\leftidx{^*}{M}{}} = 
\ifx\du\undefined
  \newlength{\du}
\fi
\setlength{\du}{10\unitlength}

\;.
\end{equation}
For a right $A$-module $(M, \rho^{M})$ and a left $A$-module $(X,\rho^{X})$, there is a notion of the tensor product  $M \tensor[A] X$ over $A$, see e.g. \cite{SchwCat}.  $M \tensor[A] X$ is an object in $\Cat{C}$ that is 
defined as the cokernel of the map $(\rho^{M} \otimes \id_{X})-(\id_{M} \otimes \rho^{X}) : M \otimes A \otimes X \rr M \otimes X$.    
When $A$ is a normalised special Frobenius algebra,   $M \tensor[A] X$ is equal to the image of the following projector $P: M \otimes X \rr M \otimes X$: 
\begin{equation}\label{equation:projector1}
  P=
\ifx\du\undefined
  \newlength{\du}
\fi
\setlength{\du}{10\unitlength}

\;.
\end{equation}
\begin{proposition}
  \label{proposition:modules-non-zero}
  Let $A$ be an algebra in a  fusion category $\Cat{C}$.  
  There is a  natural isomorphism $\Hom_{A}(M,N) \simeq \Hom_{\Cat{C}}(M \tensor[A] \leftidx{^{*}}{}{N},1)$ for $M,N \in \Mod(A)$.
\end{proposition}
\begin{proof}
  This follows from the properties of tensor product over $A$, see also \cite[Lemma 7.8.24]{EGNObook}.
 \end{proof}
\subsection{Module Categories}
\label{sec:module-categories}
 In this subsection we summarise  the main definitions and results concerning module categories, see \cite{EGNObook,Ostrik} for more details. The following definition 
is a restriction of the definition in \cite{Ostrik} to semisimple categories.
\begin{definition}
  A left $\Cat{C}$-module category $\Cat{M}$ is a semisimple $\C$-linear abelian category $\Cat{M}$,
  together  with a bifunctor $\act : \Cat{C} \times \Cat{M} \rr \Cat{M}$ and
  natural isomorphisms
  \begin{equation}
    \label{eq:structures-module-category}
    \omega_{c,d,m}: (c\otimes d) \act m \rr c\act (d \act m), \quad
    l_M: 1\act m \rr m, 
  \end{equation}
  for all $c, d \in \Cat{C}$, $m \in \Cat{M}$, such that the module constraints are fulfilled: The diagrams
  \begin{equation}
    \label{eq:diagramm}
    \begin{xy}
      \xymatrix{
     & ((c \otimes d) \otimes e) \act m \ar[ld]_{=}
        \ar[rd]^{\omega_{c\otimes d,e,m}}& \\
       (c \otimes (d \otimes e)) \act m \ar[d]^{\omega_{c,d\otimes e,m}} & & (c
        \otimes d ) \act (e \act m) \ar[d]^{\omega_{c,d, e\act m}}\\
        c \act (( d \otimes e) \act m)\ar[rr] ^{\id_c\act \omega_{d,e,m}}& & c \act (d \act (e \act m)),
      }
    \end{xy}   
  \end{equation}
  and
  \begin{equation}
    \label{eq:triangle}
    \begin{xy}
      \xymatrix{
        (c\otimes 1) \act m \ar[rr]^{ \omega_{c,1,m}}\ar[dr]^{\id_c \act m} & & c \act (1 \act  m) \ar[dl]^{\id_c \act l_m}\\
        & c \act m & 
      }
    \end{xy}
  \end{equation}
  commute for all objects $c,d,e \in \Cat{C}$ and $m \in \Cat{M}$. To emphasise that $\Cat{M}$ is a left $\Cat{C}$-module category we denote it $\CM$. 
  There is an analogous definition of a right $\Cat{C}$-module category $\MC$ with an bifunctor $\ract: \MC \otimes
  \Cat{C} \rr \MC$ satisfying analogous constraints. 
\end{definition}
For a left $\Cat{C}$-module category $\CM$, the opposite  category $\op{M}$ is a right $\Cat{C}$-module category $\MCop$ with action 
\begin{equation}
  \label{eq:from-left-to-right-module-cat}
  m\ractopp c= c^{*} \act m.
\end{equation}

\begin{definition}{\rm \cite{Ostrik}}
 Let $\CM$ and $\CN$ be $\Cat{C}$-module categories.
  \begin{definitionlist}
  \item  A $\Cat{C}$-module functor ${F}: \CM \rr\CN$ is a
    functor ${F}$ together with natural isomorphisms $f_{c,m}: F(c \act
    m ) \rr c  \act F(m)$, such that the usual pentagon and triangle diagrams commute, see \cite{Ostrik}.  
    We sometimes write $(F,f)$ for a
    module functor and call $f$ a left module constraint for $F$. 
    Module functors between right $\Cat{C}$-module categories are defined analogously.
  \item Let  $(F,f) :\CM\rr
    \CN$ and $(G,g): \CM \rr \CN$ be module functors.  A module natural transformation $ \eta: F \rr G$ is a natural
    transformation for which the diagrams
    \begin{equation}
      \label{eq:module-nat-transf}
      \begin{xy}
        \xymatrix{
          F(c \act m) \ar[r]^{\eta(c \act m)} \ar[d]^{f_{c,m}} & G(c \act m) \ar[d]^{g_{c,m}} \\
          c \act F(m) \ar[r]^{\id_c \act \eta(m)} & c \act G(m),
        }
      \end{xy}
    \end{equation}
    commute for all possible objects. The category of module functors from $\CM$ to $\CN$ and module natural transformations between them is denoted by $\Funl{\Cat{C}}{\CM,\CN}$.
  \end{definitionlist}
\end{definition}
It is easy to see that the adjoint functor of a module functor is again a module functor.  Its module functor constraint is uniquely determined by the requirement that the evaluation and coevaluation of the adjunction
are module natural transformations.
Two module categories $\CM$ and $\CN$ over $\Cat{C}$ are called equivalent if there exist module functors  $(F,f) :\CM\rr
\CN$ and $(G,g): \CN \rr \CM$ and module natural isomorphisms $F\circ G \rr \id_{\Cat{N}}$ and $G \circ F \rr \id_{\Cat{M}}$.
The 2-category of left module categories over $\Cat{C}$, module functors and module natural transformations between them is called $\ModC$.

There is an obvious notion of a   
submodule category and of a direct sum  of module categories. 
A module category is called indecomposable if it is not equivalent to a direct sum of two non trivial module categories, 
and it is called irreducible if it has no nontrivial submodule categories.
It is shown  in \cite[Lemma 1]{Ostrik} that
a module category $\Cat{M}$ over $\Cat{C}$ is indecomposable if and only if it is irreducible and that in this case there are finitely many isomorphism classes of simple objects in $\Cat{M}$.
In particular, there exists a complement  for every submodule category.

The category  of modules over a separable algebra $A \in \Cat{C}$ is a  $\Cat{C}$-module category with action given by  equation (\ref{eq:tensor-module}). It is indecomposable if and only if the algebra is indecomposable \cite[Remark 5]{Ostrik}.
The following theorem leads to the notion of Morita equivalence of fusion categories.
\begin{theorem}{\rm \cite{Mue1,ENOfus}}
  \label{theorem:Dual-is-fusion}
  Let $\CM$ be an indecomposable left $\Cat{C}$-module category. The category of $\Cat{C}$-module functors  $\Funl{\Cat{C}}{\CM,\CM}$ is a fusion category with monoidal structure given by composition of functors and 
  duality by the adjunction of module functors.
\end{theorem}
$\Funl{\Cat{C}}{\CM,\CM}$ is called the category dual to $\Cat{C}$ with respect to $\CM$. 
In particular, all module natural isomorphisms from the identity functor of an indecomposable module category to itself are multiples of the identity.

In the sequel we also require the notion of balanced functors between module categories. 
\begin{definition}{\rm \cite{Green}}
  \label{definition:Balanced-module-functor}
  Suppose $\MC$ is a right $\Cat{C}$-module category, $\CN$ a left $\Cat{C}$-module category and $\Cat{A}$ and additive category. 
  \begin{definitionlist}
  \item A functor $F: \Cat{M} \times \Cat{N} \rr \Cat{A}$ is called
    $\Cat{C}$-balanced if it is equipped with natural isomorphisms $f_{m,c,n}:F(m \ract c,n)\rr F(m,c \act n)$ for all objects $c \in \Cat{C}$, $m \in \Cat{M}$ and $n \in \Cat{N}$, such that the pentagons 
    \begin{equation}
      \label{eq:C-balanced-functor}
      \begin{xy}
        \xymatrix{
          & F(m \ract (c \otimes d),n) \ar[dl] \ar[dr]^{f_{m,c\otimes d,n}}  & \\
          F((m\ract c)\ract d,n)\ar[d]^{f_{m \ract c, d,n}} & & F(m, (c \otimes d) \act n) \ar[d]\\
          F(m \ract c, d \act n) \ar[rr]^{f_{m,c,d \act n}}& & F( m, c \act (d \act n)),
        }
      \end{xy}   
    \end{equation}
    commute for all possible objects. The unlabelled lines are the isomorphisms obtained from the module constraints of $\Cat{M}$ and $\Cat{N}$, respectively. The natural isomorphism $f$ is called balancing constraint.
  \item  Let $F,G: \Cat{M} \times \Cat{N} \rr \Cat{A}$ be two $\Cat{C}$-balanced functors with balancing constraints $f$ and $g$, respectively. A $\Cat{C}$-balanced natural transformation $\eta:F \rr G$ is a natural 
    transformation, such that the diagrams 
    \begin{equation}
      \label{eq:balancing-nat-trans}
      \begin{xy}
        \xymatrix{
          F(m \ract c,n) \ar[r]^{\eta(m \ract c, n)} \ar[d]^{f_{m,c,n}} & G(m \ract c,n) \ar[d]^{g_{m,c,n}} \\
          F(m, c\act n) \ar[r]^{\eta(m, c \act n)}  & G(m, c \act n)
        }
      \end{xy}
    \end{equation}
    commute for all possible objects. 
  \end{definitionlist}
\end{definition}
\section{Module Traces}
\label{sec:module-traces}
In this section we introduce a graphical calculus for module categories and derive a graphical  description of the algebra morphism of the inner hom objects. In the second subsection 
we introduce module traces and discuss their basic properties and some examples. 
\subsection{ Graphical Calculus for Module Categories}
\label{subsec:Inner-hom}

We extend the   graphical calculus for tensor categories (see Appendix \ref{subsection:Graphical-calculus}) to module categories. We represent objects, morphisms and the action on a module category
$\CM$ as follows.
\begin{equation}
  m \corres 
\ifx\du\undefined
  \newlength{\du}
\fi
\setlength{\du}{10\unitlength}

\;.
\end{equation}
Any module category is equivalent to a strict module category, see \cite[Thm. 1.3.8.]{GreenPhd}. 
This implies that the  graphical notation for module categories has properties  analogous to the graphical notation for tensor categories: Once  parentheses and actions of unit objects are specified 
for the incoming and outgoing objects, each diagram unambiguously represents  a  morphism in $\Cat{M}$.

We briefly summarise the definition of the inner hom object from \cite{Ostrik}.
Let $\Cat{M}$ be a left $\Cat{C}$-module category. An inner hom object $\icmM{m,n} \in \Cat{C}$ for $m,n \in \Cat{M}$ is an object in $\Cat{C}$ with  a collection of isomorphisms
\begin{equation}
  \label{eq:inner-hom}
  \alpha: \Hom_{\Cat{M}}(c \act m,n) \simeq \Hom_{\Cat{C}}(c, \icmM{m,n}), 
\end{equation}
that is natural in  $c \in \Cat{C}$.  We write $\icm{.,.}$ when the relevant module category $\Cat{M}$ is clear from the context.
Inner hom objects always exist, are unique up to a unique isomorphism and determine a bifunctor 
 $ \icmM{.,.}: \op{M} \times \Cat{M} \rr \Cat{C}$, such that the isomorphism (\ref{eq:inner-hom}) is natural in all arguments.  In the following we will speak of ``the inner hom object''.
The inner hom bifunctor is compatible with the module structure  \cite{Ostrik}:
\begin{equation}
  \label{eq:inner-hom-comp-module}
  \icm{m,c \act n} \simeq  c \otimes \icm{m,n}, \quad \text{and} \quad \icm{c \act m,n} \simeq \icm{m,n} \otimes c^{*}.
\end{equation}

we represent the inner hom object  by the following diagram: 
\begin{equation}
  \icm{m,n} \corres 
\ifx\du\undefined
  \newlength{\du}
\fi
\setlength{\du}{10\unitlength}

\;,
\end{equation}
and the isomorphism (\ref{eq:inner-hom}) reads:
\begin{equation}
  \alpha: 
\ifx\du\undefined
  \newlength{\du}
\fi
\setlength{\du}{10\unitlength}
%
\;.
\end{equation}
This can be visualised by flipping the string representing $m$ and zipping it with the $n$-string.
For a morphism $g: n \rr \tilde{n}$,  the morphism $\icm{m,g}: \icm{m,n} \rr \icm{m,\tilde{n}}$  is given  by the diagram 
\begin{equation}
\ifx\du\undefined
  \newlength{\du}
\fi
\setlength{\du}{10\unitlength}

\;.   
\end{equation}
Each morphism $h: m \rr \tilde{m}$ defines a morphism $\icm{h,n}: \icm{\tilde{m},n} \rr \icm{m,n}$ that is depicted as
\begin{equation}
\ifx\du\undefined
  \newlength{\du}
\fi
\setlength{\du}{10\unitlength}

\;.
\end{equation}
  The symbol $h^{*}$ indicates that the functor $\icm{.,.}$
is contravariant in the first argument. 
\begin{remark}
In the case of  $\Cat{C}$ considered  as a left module category over itself,
 the inner hom object of $c,d \in \Cat{C}$ is given by $\icm{c,d}=d \otimes c^{*}$.
For a morphism $h: c \rr \tilde{c}$ indeed $\icm{h,d}= \id_{d} \otimes h^{*}$. The notation $h^{*}$ therefore is consistent.
\end{remark}
The naturality of $\alpha: \Hom_{\Cat{M}}(c \act m,n) \simeq \Hom_{\Cat{C}}(c, \icm{m,n}) $ manifests itself in the graphical calculus as follows:
\begin{enumerate}
\item $\alpha$ is natural with respect to $m$:
  \begin{equation}
    \alpha: 
\ifx\du\undefined
  \newlength{\du}
\fi
\setlength{\du}{10\unitlength}

\;.
  \end{equation}
\item $\alpha$ is natural with respect to $n$:
  \begin{equation}
    \alpha: 
\ifx\du\undefined
  \newlength{\du}
\fi
\setlength{\du}{10\unitlength}
%
\;.
  \end{equation}
\item $\alpha$ is natural with respect to $c$:\
  \begin{equation}
    \alpha: 
\ifx\du\undefined
  \newlength{\du}
\fi
\setlength{\du}{10\unitlength}
%
\;.
  \end{equation}
  \end{enumerate}
\begin{lemma}
  \label{lemma:alpha-comp-tensor}
  The natural isomorphism $\alpha$ from equation (\ref{eq:inner-hom}) is compatible with the module structure. For all morphisms $\gamma : x \rr y$ in $\Cat{C}$ and 
  all $f \in \Hom_{}(c \act m, n)$,
  \begin{equation}
\ifx\du\undefined
  \newlength{\du}
\fi
\setlength{\du}{10\unitlength}

\;.
  \end{equation}
\end{lemma}

\begin{proof}
  It suffices to proof the statement for $y=x$ and $\gamma=\id_{x}$. The general case then follows directly from the naturality of $\alpha$. 
  First recall that the canonical isomorphism $\icm{m,c\act n} \simeq c \otimes \icm{m,n}$ is constructed as follows.
  Consider for $x,c \in \Cat{C}$ and $m,n \in \Cat{M}$ the square: 
  \begin{equation}
    \begin{xy}
      \xymatrix{
        \Hom_{}(x \act (c \act m),n) \ar[r]^{\simeq} \ar[d]^{\simeq}& \Hom_{}(c \act m, \leftidx{^{*}}{x}{} \act n) \ar[d]^{\alpha} \\ 
        \Hom_{}((x \otimes c)\act m,n)\ar[d]^{\alpha} &\Hom_{}(c, \icm{m,  \leftidx{^{*}}{x}{} \act n}) \ar[d]^{\simeq} \\
        \Hom_{}(x \otimes c, \icm{m,n}) \ar[r]^{\simeq}  &\Hom_{}(c,  \leftidx{^{*}}{x}{} \otimes \icm{m,n} ).
      }
    \end{xy}
  \end{equation}
  The horizontal isomorphisms are induced by the duality in $\Cat{C}$, while the unlabelled vertical isomorphism on the right is the
natural isomorphism $ \icm{m, \leftidx{^{*}}{x}{} \act n} \simeq  \leftidx{^{*}}{x}{} \otimes \icm{m,n}$ from equation (\ref{eq:inner-hom-comp-module}).
  This isomorphism is defined by the requirement that  the square commutes. As it is constructed from 
  natural isomorphisms which we suppress in the graphical notation,  we will suppress this isomorphism as well in the sequel.
  It follows from   the commutativity of the previous diagram, that  the diagram
  \begin{equation}
    \begin{xy}
      \xymatrix{ \Hom_{}(x \act (c \act m), x \act n) \ar[r]^{\simeq} \ar[d]^{\alpha} & \Hom_{}(x^{*} \otimes x  \otimes c) \act m,  n) \ar[d]^{\alpha} \\
        \Hom_{}(x \otimes c , x \otimes \icm{m,n}) & \Hom_{} ((x^{*} \otimes x  \otimes c, \icm{m,n}) \ar[l]^{\simeq}
      }
    \end{xy}
  \end{equation}
  commutes. If we choose $\id_{x} \otimes f \in  \Hom_{}(x \act (c \act m), x \act n) $ with $f \in \Hom_{}(c \act m,n)$ in the left upper space, 
the commutativity of the diagram implies $\alpha(\id_{x} \otimes f)=\id_{x} \otimes \alpha(f)$.
\end{proof}

A useful property of the inner hom is that  $\icm{m,m}$ has a canonical structure of an algebra in $\Cat{C}$ \cite{Ostrik}. Next we present a graphical definition of this structure. 
The internal evaluation morphism $\ev{n,m}: \icm{n,m} \act n \rr m$ (see \cite[Sec. 3.2.]{Ostrik}) is 
defined by:
\begin{equation}
  \ev{n,m} =\alpha^{-1}(\id_{\icm{n,m}}) \quad \corres \quad  
\ifx\du\undefined
  \newlength{\du}
\fi
\setlength{\du}{10\unitlength}

\;.
\end{equation}
This notation is compatible with the notation for $\alpha$ since by flipping the $n$-string we obtain the identity string $\icm{n,m}$.
The internal multiplication $\mu_{m,n,k}: \icm{n,k} \otimes \icm{m,n} \rr \icm{m,k}$ and the internal unit $\eta_{m}: 1 \rr \icm{m,m}$ are given by
\begin{equation}
  \mu_{m,n,k} =\alpha \left( 
\ifx\du\undefined
  \newlength{\du}
\fi
\setlength{\du}{10\unitlength}

\;.
\end{equation}

\begin{lemma}
  \label{lemma:technical-inner-hom}
  For all morphisms $f \in \Hom_{}(c \act m, \icm{n,k}\act n)$, 
  \begin{equation}
    \alpha( \ev{n,k} \circ f )= \mu_{m,n,k} \circ \alpha(f), \;\text{i.e.}
  \end{equation}
  \begin{equation}
\ifx\du\undefined
  \newlength{\du}
\fi
\setlength{\du}{10\unitlength}

\;,
  \end{equation}
  for all $g \in \Hom_{}( d \act l, s)$
  follows from applying $\alpha$ to both sides and using the naturality of $\alpha$. Applying this identity 
  to $f$ with $s= \icm{n,k}\act n $ yields
  \begin{equation}
\ifx\du\undefined
  \newlength{\du}
\fi
\setlength{\du}{10\unitlength}

\;.
  \end{equation}  
  Applying $\alpha$ to  the right hand side of this equation and  using its naturality proves the claim.
\end{proof}
 The following theorem plays an important role in the 
theory of fusion categories since it combines the theory of module categories with the theory of algebras. 
\begin{theorem}{\rm \cite{Ostrik}}
  \label{theorem:inner-hom-prop}
For all non-zero objects $m,n$ in a $\Cat{C}$-module category $\CM$,
  $\icm{m,m}$ is an algebra object in $\Cat{C}$ and
  $\icm{m,n}$ is a right $\icm{m,m}$-module.
  The functor $\Cat{M} \ni n \mapsto \icm{m,n} \in  \Mod(\icm{m,m})$ yields an equivalence of $\Cat{C}$-module categories
  provided $\CM$ is indecomposable.
\end{theorem}
We will revisit parts of the proof of this statement with the graphical calculus. 
\begin{proposition}
  \label{proposition:inner-hom-algebra-obj}
  \begin{propositionlist}
  \item \label{item:internal-ev-module-morph} 
    The internal evaluation morphism is a module morphism:
    \begin{equation}
\ifx\du\undefined
  \newlength{\du}
\fi
\setlength{\du}{10\unitlength}

\;.
    \end{equation}
  \item The internal multiplication  is associative: 
    \begin{equation}\label{equation:multiplication-inner-hom}
\ifx\du\undefined
  \newlength{\du}
\fi
\setlength{\du}{10\unitlength}
%
\;.
    \end{equation}
\item For all non-zero  $ m \in \Cat{M}$, $\icm{m,m}$ is canonically an algebra object.
  \end{propositionlist}
\end{proposition}
\begin{proof}
  The first relation follows from applying $\alpha$ to both diagrams. 
  Both diagrams obtained in this way represent the multiplication morphism. Since $\alpha$ is an isomorphism, the preimages have to agree as well.

  To show the second part, first note that the expression on the left hand side of equation (\ref{equation:multiplication-inner-hom}) is $\alpha$ applied to 
  \begin{equation}
\ifx\du\undefined
  \newlength{\du}
\fi
\setlength{\du}{10\unitlength}

\;.
  \end{equation}
  Now apply $\alpha$ to the diagram on the right. In the upper part of the diagram this results in the morphism $\id_{\icm{l,k}} \otimes \mu_{m,n,l}$ due to Lemma \ref{lemma:alpha-comp-tensor}.  
  With Lemma \ref{lemma:technical-inner-hom} we conclude that $\alpha$ applied to this diagram yields the right hand side of equation (\ref{equation:multiplication-inner-hom}).
  The statement follows since $\alpha$ is an isomorphism.
To show the last part we only have to prove the compatibility of the internal multiplication and the internal unit. This is a direct computation in the diagrammatic calculus.  
\end{proof}

\subsection{Module Traces on Module Categories over Pivotal Fusion Categories}

We are now ready to define the notion of a  module trace.
 As an example  we discuss module categories over $G$-graded vector spaces. This illustrates that the existence of a module 
trace on a given module category distinguishes different pivotal structures.

For each module category $\Cat{M}$ over a pivotal category $\Cat{C}$ there is a linear map  
\begin{equation}
  \label{eq:trace-C-on-M}
  \tr^{\Cat{C}}_{c,m}: \End_{\Cat{M}}(c \act m) \rr \End_{\Cat{M}}(m), \quad f \mto  (\evp{c} \act \id_{m}) \circ (a_{^{*}c} \act f)  \circ(\coev{c} \act \id_{m}),
\end{equation}
which  we call partial  trace. Whenever this is unambiguous  we omit the labels of $\tr^{\Cat{C}}$.
The graphical representation of this map is 
\begin{equation}
  \tr^{\Cat{C}} \left( 
\ifx\du\undefined
  \newlength{\du}
\fi
\setlength{\du}{10\unitlength}
\;.
\end{equation}
As a direct consequence of the definition of a module functor we obtain:
\begin{lemma}\label{lemma:partial-trace-module-funct}
  Let $F: \CM \rr \CN$ be a $\Cat{C}$-module functor. For all $f \in \End_{\Cat{M}}(c \act m)$, $\tr^{\Cat{C}}(F(f))= F(\tr^{\Cat{C}}(f))$.
\end{lemma}
With the map $\tr^{\Cat{C}}$ we can define module traces.
\begin{definition}\label{definition:Module-trace}
  Let $\Cat{M}$ be a module category over a pivotal fusion category $\Cat{C}$. A  trace $\Theta$ on $\Cat{M}$ is a collection of linear maps 
  \begin{equation}
    \label{eq:theta}
    \Theta_{m}: \End_{\Cat{M}} (m)\rr \C \quad \text{for all} \quad m \in \Cat{M},
  \end{equation}
  such that the following properties are satisfied:
  \begin{definitionlist}
  \item 
    \label{item:Theta-symmetric} 
    $\Theta$ is symmetric: for all $f \in \Hom_{\Cat{M}}(m,n)$ and $g\in \Hom_{\Cat{M}}(n,m)$,
    \begin{equation}
      \label{eq:symmetric}
      \Theta_{m}(g \circ f) = \Theta_{n}(f \circ g).
    \end{equation}
  \item 
    \label{item:non-degenerate}$\Theta$ is non-degenerate: the pairing 
    \begin{equation}
      \label{eq:non-deg}
      \Hom_{\Cat{M}}(m,n) \times \Hom_{\Cat{M}}(n,m) \rr \C, \quad (f, g) \mto \Theta_{m}(g \circ f)
    \end{equation}
    is non-degenerate for all $m,n \in \Cat{M}$.
  \end{definitionlist}
  If furthermore 
  \begin{definitionlist}
    \addtocounter{enumi}{2}
  \item
    \label{item:C-compatible} 
    $\Theta$ is $\Cat{C}$-compatible: for all $c \in \Cat{C}$, $m \in \Cat{M}$, 
    \begin{equation}
      \label{eq:Theta-C-module}
      \Theta_{c\act m}= \Theta_{m} \circ \tr^{\Cat{C}}, 
    \end{equation}
  \end{definitionlist}
  then $\Theta$ is called a $\Cat{C}$-module trace or module trace if the category $\Cat{C}$ is clear from the context. 
  We sometimes write $(\MC, \Theta^{\Cat{M}})$ for a module category with module trace.
\end{definition}
The notion of a trace on a linear category is well-known and  a category with a trace is  also called a Calabi-Yau category, see e.g. \cite[Sec. 2]{Cost}.
\begin{remark}\label{remark:Theta}
  \begin{remarklist}
  \item \label{item:example-pivotal}
    The  notion of a module trace  is a generalisation of the trace on  a pivotal fusion category. Indeed, consider $\Cat{C}$ as a left module category over itself.
    The left trace $\tr^{L}: \End_{\Cat{C}}(c) \rr \C$ induces a canonical  module trace on $\Cat{C}$.
    The  left trace is symmetric and the  compatibility of the duality with the tensor product yields equation (\ref{eq:Theta-C-module}).
    For the non-degeneracy, note that the argument in the proof of  \cite[Lemma \Rmnum{2}.4.2.3]{Tur} can be extended  to the case of pivotal fusion categories, 
see also Lemma \ref{lemma:general-trace-ses}. The semisimplicity of $\Cat{C}$ is crucial at this point and for this reason 
     we restrict attention to fusion categories in this work, see \cite{Pivthree}  for a generalization to  more general tensor categories.

  \item For any given  trace $\Theta$ on $\Cat{M}$ and non-zero number $z \in \C$ the linear maps $z \cdot \Theta_{m}$ define 
    another  trace denoted  $z \cdot \Theta$. If $\Theta$ is  a module trace then  $z \cdot \Theta$ is again a module trace. We will show in Section \ref{subsection:Characterization-module-trace}
 that module traces are unique up to such  rescalings. 
  \end{remarklist}
\end{remark}
We introduce a graphical notation for module traces:
\begin{equation}
  \Theta_{m}(f) \quad \corres \quad
\ifx\du\undefined
  \newlength{\du}
\fi
\setlength{\du}{10\unitlength}
\;.
\end{equation}
The symmetry  and  $\Cat{C}$-compatibility of $\Theta$ then read
\begin{equation}
\ifx\du\undefined
  \newlength{\du}
\fi
\setlength{\du}{10\unitlength}
%
 \;. 
\end{equation}
Given a  trace $\Theta$ on a category $\Cat{M}$, we define the dimensions of objects $m \in \Cat{M}$ with respect to $\Theta$ as 
\begin{equation}
  \label{eq:dimension}
  \dim^{\Theta}(m)=\Theta_{m}(\id_{m}).
\end{equation}
The dimensions depend only on the isomorphism classes of objects:
\begin{lemma}
  \label{lemma:first-implications}
  \begin{lemmalist}
  \item \label{item:dim-well-defined} If two objects $m,n \in \Cat{M}$ are isomorphic then $\dim^{\Theta}(m)=\dim^{\Theta}(n)$.
  \item \label{item:cptl-direct-sum} 
    $\Theta$ is compatible with direct sums. For all $f \in \End_{\Cat{M}}(m)$, $g \in \End_{\Cat{M}}(n)$, we have
    $\Theta_{m\oplus n}(f \oplus g)=\Theta_{m}(f)+ \Theta_{n}(g)$. 
    In particular,   $\dim^{\Theta}(m \oplus n)= \dim^{\Theta}(m) + \dim^{\Theta}(n)$.
  \item \label{item:dimensions-after-action}
    \begin{equation}
      \label{eq:action-dimension}
      \dim^{\Theta}(c \act m)=\dim^{\Cat{C}}(c) \cdot \dim^{\Theta}(m).
    \end{equation}
  \end{lemmalist}
\end{lemma}
\begin{proof}
  For the first part  choose an isomorphism $f:m \rr n$. The symmetry of $\Theta$ implies
  \begin{equation}
    \dim^{\Theta}(m)=\Theta_{m}(\id_{m})=\Theta_{m}(f^{-1} \circ f)=\Theta_{n}(f \circ f^{-1})= \Theta_{n}(\id_{n})= \dim^{\Theta}(n).
  \end{equation}
The second part follows directly from the linearity of $\Theta$.
  The third part is a consequence of the  $\Cat{C}$-compatibility of $\Theta$.
\end{proof}
\paragraph{Direct Sums and Equivalences of Module Categories with Module Trace}
We show that the notion of a module trace is well-behaved with respect to decomposition of module categories
and investigate the structure of the module categories with $\Cat{C}$-module trace in the 2-category $\ModC$ of $\Cat{C}$-module categories, module functors and module natural transformations. 
\begin{definition}
  Let $\ModCt$ be the full sub 2-category of $\ModC$ which has $\Cat{C}$-module categories $(\MC, \Theta^{\Cat{M}})$ endowed with a $ \Cat{C}$-module trace $\Theta^{\Cat{M}}$ as objects. 
  A module functor $F: \CM \rr \CN$ is called an 
  isometric module functor if $\Theta^{\Cat{N}}(F(f))= \Theta^{\Cat{M}}(f)$ for all $f \in \End_{\Cat{M}}(m)$ and all $m \in \Cat{M}$. 
  Two module categories in $\ModCt$ are called isometrically equivalent 
  if there exists an equivalence of module categories consisting of isometric module functors between them.
\end{definition}
Note that an isometric module functor is faithful due to the non-degeneracy of the module traces.
The subcategory  $\ModCt$ is well-behaved in the following sense. 
\begin{proposition}
  \label{proposition:traces-well-behaved-in-Mod}
  \begin{propositionlist}
  \item 
    Let $(\CM, \Theta^{\Cat{M}})$ be a an object in $\ModCt$ and let $\CN$ be a module category with an equivalence  $F: \CN \rr \CM$ of module categories. 
    Then there exists a  $\Cat{C}$-module trace on $\CN$ such that $F$ is an isometric equivalence.
  \item The direct sum of two  module categories with module traces possesses a canonical module trace. 
  \item A submodule category of a module category with module trace 
    inherits a canonical module trace. 
  \item Each object in $\ModCt$ is isometrically equivalent to a finite direct sum of indecomposable objects.
  \end{propositionlist}
\end{proposition}
\begin{proof}
To show the first part, define the linear maps  $\Theta^{\Cat{N}}(f)=\Theta^{\Cat{M}}(F(f))$ for all $f \in \End_{\Cat{N}}(n)$.  Lemma \ref{lemma:partial-trace-module-funct} implies that this defines a module 
trace for $\Cat{N}$ and that  $F$ is isometric by construction.
  For the second part consider an object $m \oplus n \in \CM \oplus \CN$. Since  $\End_{\Cat{M} \oplus \Cat{N}}(m \oplus n)= \End_{\Cat{M}}(m) \oplus \End_{\Cat{N}}(n)$, 
  we can define a linear map $(\Theta^{\Cat{M}} \oplus \Theta^{\Cat{N}})_{m \oplus n}: \End_{\Cat{M} \oplus \Cat{N}} (m \oplus n)\rr \C $ as the sum $\Theta^{\Cat{M}}_{m} \oplus \Theta^{\Cat{N}}_{n}$. 
  It is easy to see that this defines a $\Cat{C}$-module trace. 
 Now consider a submodule category of a module category with module trace. As we can choose a complement of the submodule category, the restriction of a
  module trace to a submodule category is non-degenerate and hence a module trace.
  The last statement is a consequence of the first and second statement.
\end{proof}

\paragraph{Examples\\}
We denote by $\Vect$  the fusion category of finite dimensional $\C$-vector spaces. It has  a unique pivotal structure.
 A semisimple abelian category over $\C$ is  a  module category over $\Vect$ 
with module structure  $V \tensor[\C]m$ defined by   $ V \tensor[\C] \Hom_{\Cat{M}}(m,n) \simeq \Hom_{\Cat{M}}(V \tensor[\C]m,n)$ for  $V \in \Vect$ and $m,n \in \Cat{M}$.
\begin{lemma}
  \label{lemma:Vect-mdoule-trace}
  A trace on a semisimple category $\Cat{M}$ is also a $\Vect$-module trace on $\Cat{M}$.
\end{lemma}
\begin{proof}
  We show that condition (\ref{eq:Theta-C-module}) is satisfied. As $\End_{\Cat{M}}(V \tensor[\C]m) \simeq \End(V) \otimes \End(m)$, it is sufficient to show that 
  \begin{equation}
    \label{eq:technichal}
    \Theta_{V \tensor[\C] m}(\alpha \tensor[\C] f)=\tr(\alpha) \Theta_{m}(f)
  \end{equation}
for all $V \in \Vect$,
  $\alpha \in \End(V)$ and $f \in \End(m)$. Here $\tr$ is the usual trace on $\Vect$ that coincides with the left trace on  $\Vect$ considered as a fusion 
  category. 
 Equation (\ref{eq:technichal}) follows from a  direct calculation in  a basis for $V$.
\end{proof}

\begin{example}
  \label{example:finite-G}
  Let $G$ be a finite group and $\omega \in \Chain^{3}(G, \C^{\times})$  a normalised cocycle.
  This data defines a fusion category $\Vect_{G}^{\omega}$ with simple objects labelled by elements of $G$, see \cite{EGNObook} and \cite{Ostrik} for more details. 
  The pivotal structures on $G$ are in bijection with the characters  $ \kappa \in \Hom(G, \C^{\times})$. Indecomposable module categories $M(H, \Psi)$ 
  over  $\Vect_{G}^{\omega}$ are obtained from subgroups $H \subset G$ with $\omega |_{H}=1$ and 
  cocycles $\Psi \in \Chain^{2}(H, \C^{\times})$. The simple objects of $M(H, \Psi)$ are labelled by elements in the right cosets $[ g ] \in H \backslash G$. The action of a simple object $x \in \Vect_{G}^{\omega}$
  is given by $x \act [g]=[xg]$, with module constraint twisted by $\Psi$. 

 A module category  $M(H, \Psi)$ over the pivotal fusion category $(\Vect_{G}^{\omega}, \kappa)$ possesses a module trace if and only if $\kappa|_{H}=1$. 
  This can be seen as follows: Suppose $\Theta$ is a module trace on  $M(H, \Psi)$ normalised by $\Theta([e])=1$. 
  Then equation (\ref{eq:action-dimension}) implies  $\Theta([gx])=\kappa(g)\cdot \Theta([x])$, in particular $\Theta([g])=\kappa(g)$. So $\kappa$ is well-defined on  $ H \backslash G$, which is the case if and only if 
  $\kappa|_{H}=1$. Conversely, if  $\kappa|_{H}=1$ it is easy to see that $\kappa$ yields a module trace for  $M(H, \Psi)$.  In particular, there exists a module trace for all module categories over $\Vect_{G}$ when the 
  pivotal structure $\kappa \equiv 1$ is chosen. 
\end{example}

\begin{example}
  Let $\Cat{C}$ be a fusion category. Recall the construction of a pivotal fusion category $\widetilde{\Cat{C}}$ from \cite[Remark 3.1]{ENOfus}: 
  The simple objects of $\widetilde{\Cat{C}}$ are pairs $(c,f_{c})$, where $c \in \Cat{C}$ is a simple object and $f_{c}: c \rr c^{**}$ is an isomorphism such that $f_{c}^{**}f_{c}=g_{c}$, where $g$ is the canonical 
  monoidal natural isomorphism $\id_{\Cat{C}} \rr (.)^{****}$ defined in \cite{ENOfus}. With  $(c,f_{c}) \in \widetilde{\Cat{C}} $, also  $(c,-f_{c}) \in \widetilde{\Cat{C}}$.
  $\widetilde{\Cat{C}}$ has a canonical pivotal structure such that $\dim^{\widetilde{\Cat{C}}}(c,f_{c})=\evp{c^{**}} \circ (f_{c} \otimes \id_{c^{*}})\circ \coev{c}=:   \tr(f_{c})$. 
  The monoidal structure of  $\widetilde{\Cat{C}}$ is induced by the monoidal structure of $\Cat{C}$ and 
  the forgetful functor $U :\widetilde{\Cat{C}} \rr \Cat{C}$ is a monoidal functor. Hence $\Cat{C}$ is a left  $\widetilde{\Cat{C}}$-module category. This module category does not admit a module trace when  $\widetilde{\Cat{C}}$ 
  is equipped with the canonical pivotal structure: 
  Assume that $\Theta$ is a module trace and let $(c,f_{c}) \in \widetilde{\Cat{C}}$ and $d \in \Cat{C}$ 
be simple objects. Then
  \begin{equation}
    \begin{split}
      \tr(f_{c})\dim^{\Theta}(d)&=\dim^{\Theta}((c,f_{c})\act d)=\dim^{\Theta}(c \otimes d) \\
      &=\dim^{\Theta}((c,-f_{c})\act d)=-\tr(f_{c}) \dim^{\Theta}(d),
    \end{split}
  \end{equation}
 which is a contradiction, since as we will explain in Section \ref{sec:existence-problem-as}, $ \dim^{\Theta}(d) \neq 0$.

However, a pivotal structure $a$ for $\Cat{C}$ induces  a different pivotal structure for  $\widetilde{\Cat{C}}$ with quantum dimensions
  $\dim^{(\widetilde{\Cat{C}},a)}(c,f_{c})=\dim^{\Cat{C}}(c)$  and it is easy to see that the left trace with respect to $a$ defines a $\widetilde{\Cat{C}}$-module trace for the module category $\Cat{C}$. 
\end{example}

These examples motivate the following definition. 
\begin{definition}
  Let $\Cat{C}$ be a fusion category with pivotal structure $a$ and $\CM$ a module category.  The pair $(a, \CM)$ is called matched if there exists a $\Cat{C}$-module trace on 
  $\CM$. A pivotal structure for $\Cat{C}$ that is matched with all module categories is called flexible. 
\end{definition}
In Proposition \ref{proposition:pseudo-unitary} we will show that a pseudo-unitary $\Cat{C}$ admits a flexible pivotal structure that is also spherical. 
It has been conjectured in \cite{ENOfus} that all fusion categories admit a pivotal structure. The theory of module traces raises the following refinements of this question.

\begin{enumerate}[(1)]
\item Given  a fusion category  $\Cat{C}$ and an indecomposable module category $\CM$, is there  a pivotal structure $a$ on $\Cat{C}$, such that the pair $(a, \CM)$ is matched? 
\item  Does every (modular) fusion categories  exhibit  a   flexible pivotal structure   and is it unique?
\item Is every flexible pivotal structure  spherical?
\end{enumerate}

After the publication of this article, we realized that there are rather immediate answers to 
 the questions  and we are grateful to Pavel Etingof for sharing his answers to the questions  and his approach 
to question (1),  which  reduces the questions above to the conjecture \cite{EGNObook} that also every multi-fusion category (with possibly non-simple unit) has a pivotal structure: 

\begin{enumerate}[(1)] 
\item (P. Etingof \cite{EEigenval}) We can consider the multi-fusion category $\Cat{D}=\End_{\Cat{C}}(\CC \oplus \CM)$. It is immediate that a pivotal structure on $\Cat{D}$ gives a pivotal structure 
on $\Cat{C}$ and via the left trace on $\Cat{D}$ restricted to $\Cat{M}$ a module trace on $\CM$.
\item Every fusion category $\Cat{C}$ has finitely many indecomposable module categories 
$\CMi$, $i \in \{1, \ldots, n\}$ up to equivalence \cite{Ostrik}. 
Thus we can consider the multi-fusion category 
\begin{equation}
  \Cat{D}= \End_{\Cat{C}}(\CC \oplus \CMone \oplus \ldots \oplus \CMn ),
  \end{equation}
and a pivotal structure on $\Cat{D}$ gives a pivotal structure on $\Cat{C}$ that is matched to all 
$\CMi$, thus it is flexible. 

Let $p$ be prime. Any pivotal structure on $\Cat{C}=\Vect_{\Z_{p}}^{\omega}$ with non-trivial 3-cocycle $\omega$ is flexible: It follows from the Example \ref{example:finite-G} that the only  indecomposable 
module category over $\Cat{C}$ is $\Cat{C}$ itself which has a module trace for every pivotal structure on $\Cat{C}$. Thus the uniqueness part of the question (2) above fails, as e.g. 
$\Vect_{\Z_{2}}^{\omega}$ has already two flexible pivotal structures and is even modular. 

\item In the case $\Vect_{\Z_{3}}^{\omega}$  with non-trivial 3-cocycle $\omega$   there exist non-spherical pivotal structures which are flexible by the answer (2) above, so the answer to  question (3) is negative. 
\end{enumerate}

\section{Module Traces and  Dual $\Hom_{}$-Spaces}
\label{sec:module-traces-dual}
\subsection{Uniqueness of Module Traces}
\label{subsection:Characterization-module-trace}
In this subsection we show that  module traces are unique up to scaling.
First we examine traces on abelian categories and give an equivalent characterisation of traces in terms of certain natural isomorphisms. In the next step we show that analogous results hold for module traces.

 We denote by  $V^{*}$ the dual vector space of a $\C$-vector space $V$. 
\begin{proposition}\label{proposition:equivalences-for-trace}
  Let $\Cat{M}$ be an additive category enriched over $\Vect$. The following structures on $\Cat{M}$ are in one-to-one correspondence.
  \begin{propositionlist}
  \item A trace on $\Cat{M}$.
  \item A natural isomorphism $\eta: \Hom_{\Cat{M}}(m,n) \rr  \Hom_{\Cat{M}}(n,m)^{*}$.
  \end{propositionlist}
\end{proposition}
\begin{proof}
  Let $\Theta$ be a trace on $\Cat{M}$. The non-degenerate pairing $\Hom_{\Cat{M}}(m,n) \times
  \Hom_{\Cat{M}}(n,m)\rr \C$  defines isomorphisms $\eta_{m,n} :\Hom_{\Cat{M}}(m,n) \simeq 
  \Hom_{\Cat{M}}(n,m)^{*}$. We have to show that these isomorphisms are natural, i.e.~that for $\chi: n \rr \tilde{n}$ the diagram
  \begin{equation}
    \label{eq:natural-of-pairing}
    \begin{xy}
      \xymatrix{
        \Hom_{}(m,n) \ar[r]^{\eta_{m,n}} \ar[d]^{\Hom(m,\chi)} & \Hom_{}(n,m)^{*} \ar[d]^{\Hom(\chi,m)^{*}}\\
        \Hom_{}(m,\tilde{n}) \ar[r]^{\eta_{m,\tilde{n}}} & \Hom(\tilde{n},m)^{*}
      }
    \end{xy}
  \end{equation}
  commutes.    Let $f \in \Hom_{}(m,n)$ and $g \in \Hom_{}(\tilde{n},m)$. $\Hom(m,\chi)$ is the linear map that sends $f$ to $ \chi \circ f $. 
 In the following we denote this map by $\chi$. We compute
  \begin{equation}
    \begin{split}
      (\eta_{m,\tilde{n}}\circ \Hom(m,\chi))(f)(g)&=\Theta_{m}(g \circ (\chi \circ f))\\
      &=\Theta_{m}((g \circ \chi) \circ f)= ((\Hom(\chi,m)^{*} \circ \eta_{m,n})(f)(g).
    \end{split}
  \end{equation}
  This shows the commutativity of the diagram (\ref{eq:natural-of-pairing}).
  The proof for naturality in $m$ is analogous. 

  On the other hand, a natural isomorphism $\eta_{m,n}: \Hom_{\Cat{M}}(m,n) \rr  \Hom_{\Cat{M}}(n,m)^{*}$ induces a linear map $\Theta_{m}: \Hom_{\Cat{M}}(m,m)\rr \C $ by
  $\Theta_{m}(f)= \eta_{m,m}(\id_{m})(f)$. For $\alpha \in \Hom_{}(m,n)$ and $\beta \in \Hom_{}(n,m)$, the naturality of $\eta$ implies 
  \begin{equation}
    \Theta_{m}(\beta \circ \alpha)= \eta_{m,m}(\id_{m})(\beta \circ \alpha)= \eta_{m,n}(\alpha)(\beta)=\eta_{n,n}(\alpha \circ \beta)= \Theta_{n}(\alpha \circ \beta).
  \end{equation}
  This proves the symmetry of $\Theta$ and, as the  map $\eta_{m,n}$ is an isomorphism, also the non-degeneracy.
\end{proof}

We will now  generalise this proposition to $\Cat{C}$-module traces.
Let $\CM$ be a $\Cat{C}$-left module category. The functors
\begin{equation}
  \label{eq:functors}
  \begin{split}
    &  \op{M} \times \Cat{M} \rr \Vect, \quad {m} \times n \mapsto \Hom_{\Cat{M}}({m},n) \quad \text{and} \\
    &  \op{M} \times \Cat{M} \rr \Vect, \quad {m} \times n \mapsto \Hom_{\Cat{M}}(n,{m})^{*}
  \end{split}
\end{equation}
are canonically $\Cat{C}$-balanced (see Definition \ref{definition:Balanced-module-functor}). The balancing constraint for the first functor is the natural isomorphism 
\begin{equation}
  \label{eq:balancing-varphi1}
  \Hom_{\Cat{M}}(m, c \act n) \simeq \Hom_{\Cat{M}}(c^{*} \act m,n)= \Hom_{\Cat{M}}(m \ractopp c,n),
\end{equation}
that is available in any tensor category. In contrast,  the balancing constraint for the second functor,
\begin{equation}
  \label{eq:balancing-varphi2}
  \begin{split}
    & \Hom_{\Cat{M}}(c \act n, m)^{*}= \Hom_{\Cat{M}}( n, \leftidx{^{*}}{c}{} \act m)^* \\
    &\simeq \Hom_{\Cat{M}}( n,{m} \ractopp (\leftidx{^{**}}{c}{}))^* \simeq \Hom_{\Cat{M}}( n,{m} \ractopp c)^*, 
  \end{split}
\end{equation}
  involves the pivotal structure of $\Cat{C}$ in the last isomorphism.

\begin{theorem}
  \label{theorem:equivalences-for-C-trace}
  Let $\CM$ be a left module category over a pivotal fusion category $\Cat{C}$. The following  structures on $\Cat{M}$ are in canonical one-to-one correspondence.
  \begin{theoremlist}
  \item \label{item:module-trace} A $\Cat{C}$-module trace on $\Cat{M}$.
  \item \label{item:C-Serre}A $\Cat{C}$-balanced natural isomorphism $\eta: \Hom_{\Cat{M}}(m,n) \rr  \Hom_{\Cat{M}}(n,m)^{*}$.
  \end{theoremlist}
\end{theorem}

\begin{proof}
  We have to show that the isomorphisms $\eta_{m,n}:\Hom_{\Cat{M}}(m,n) \rr  \Hom_{\Cat{M}}(n,m)^{*}$ from Proposition 
  \ref{proposition:equivalences-for-trace} are $\Cat{C}$-balanced if and only if $\Theta$ is $\Cat{C}$-compatible. 
  Consider morphisms  $f \in \Hom_{}(m,c \act n)$ and $g \in \Hom_{}(c \act n, m)$. Denote by 
  $\hat{f} \in \Hom_{}(c^{*}\act m,n)$ and $\hat{g} \in \Hom_{}(n,c^{*} \act m)$ the morphisms obtained from $f$ and $g$ under the balancing isomorphisms (\ref{eq:balancing-varphi1}) and (\ref{eq:balancing-varphi2}), 
respectively. A direct computation shows that the $\Cat{C}$-balancing property of $\eta_{m,n}$ is 
  equivalent to the condition 
  \begin{equation}
    \label{eq:balancing-theta-equi}
    \Theta_{m}(g \circ f)=\Theta_{c^{*}\act m}(\hat{g} \circ \hat{f}),
  \end{equation}
 for all possible $f$ and $g$. Due to the symmetry of $\Theta$, 
  $\Theta_{c^{*}\act m}(\hat{g} \circ \hat{f})= \Theta_{n}(\tr^{\Cat{C}}(f \circ g))$,  and   we conclude that equation (\ref{eq:balancing-theta-equi}) is equivalent to the 
$\Cat{C}$-compatibility of $\Theta$. Thus the statement is proven.
\end{proof}
This implies in particular that for each pivotal fusion category $\Cat{C}$ there is a natural $\Cat{C}$-balanced isomorphism 
\begin{equation}
  \label{eq:eta-C}
  \eta^{\Cat{C}}:\Hom_{}(x,y)\rr \Hom_{}(y,x)^{*},
\end{equation}
induced by the left trace. 

In the sequel we will need the following  extension of the usual Yoneda lemma.
\begin{lemma}
  \label{lemma:basics-module-functor}
  Let $F, G: \CM \rr \CN$ be module functors. The set of $\Cat{C}$-module natural transformations $F \rr G$ is in canonical bijection with the set of $\Cat{C}$-balanced natural transformations of the two $\Cat{C}$-balanced 
  functors: 
  \begin{equation}
    \label{eq:two-Functors-From-mdoule-functors}
    \begin{split}
      \NCop \times \CM \ni & n \times m \mapsto \Hom_{\Cat{N}}(n,F(m))\in \Vect, \; \text{and} \\
      \NCop \times \CM \ni & n \times m \mapsto \Hom_{\Cat{N}}(n,G(m)) \in \Vect. \\
    \end{split}
  \end{equation}
A  $\Cat{C}$-balanced natural transformation $\hat{\eta}:\Hom_{\Cat{N}}(n,F(m)) \rr \Hom_{\Cat{N}}(n,G(m))$ is mapped to the unique  $\Cat{C}$-module natural transformation $ \eta:F \rr G$ with
  $\hat{\eta}(f)= \eta(m) \circ f$ for all $f \in \Hom_{\Cat{N}}(n,F(m))$.
For three module functors $F, G, K: \CM \rr \CN$, the $\Cat{C}$-module natural transformation $F \rr K$ corresponding to a composition $\Hom_{}(n,F(m)) \rr \Hom_{}(n, G(m)) \rr \Hom_{}(n, K(m))$ of $\Cat{C}$-balanced natural 
isomorphisms is equal to the composition of the corresponding $\Cat{C}$-module natural transformations.
\end{lemma}
\begin{proof}
The usual Yoneda lemma shows that a   transformation  $\hat{\eta}: \Hom_{\Cat{N}}(n,F(m)) \rr \Hom_{\Cat{N}}(n,G(m))$ that is natural in both arguments can be identified 
 with a natural transformation $\eta: F \rr G$. Consider the following diagram.
  \begin{equation}
    \label{equation:Yoneda-extension}
    \begin{xy}
      \xymatrix{
        \Hom_{}(n,F ( c\act m)) \ar[r]^{\eta(c\act m)} \ar[d]^{\simeq} & \Hom_{}(n, G(c \act m))\ar[d]^{\simeq}\\
        \Hom_{}(n, c \act F(m)) \ar[r]^{{c \act \eta(m)}} \ar[d]^{\simeq}& \Hom_{}(n, c \act G(m))\ar[d]^{\simeq}\\
        \Hom_{}(c^{*} \act n, F(m)) \ar[r]^{\eta(m)} & \Hom_{}(c^{*}\act n, G(m)).
      }
    \end{xy}
  \end{equation}
  The vertical isomorphisms provide the $\Cat{C}$-balancing structure
of the functor $\Hom_{}(n,F(m))$. It is easy to see that these isomorphisms 
  satisfy the pentagon constraint. An analogous consideration holds for $\Hom_{}(n,G(m))$.

  The lower rectangle in (\ref{equation:Yoneda-extension}) commutes for any natural transformation $\eta$. The outer diagram commutes if and only if the upper rectangle commutes. 
  The former commutes if and only if $\eta$ is a $\Cat{C}$-module natural transformation, while  commutativity of the latter is precisely the condition on $\eta$ to define 
a $\Cat{C}$-balanced natural isomorphism  $\Hom_{\Cat{N}}(n,F(m)) \rr \Hom_{\Cat{N}}(n,G(m))$.
The statement about the composition follows directly from the corresponding property of the  Yoneda lemma.
\end{proof}
The next result shows that  module traces are essentially unique. Consequently the existence of a module trace is
  a property of a module category over a pivotal fusion category  rather than a structure on a module category.

\begin{proposition}
  \label{proposition:Uniqueness}
  Let $(\CM,\Theta)$ be an indecomposable module category over $\Cat{C}$ with module trace. For any other module trace $\widetilde{\Theta}$ on $\CM$ there is a $z \in \C^{\times}$ such that $\widetilde{\Theta}=z \cdot \Theta$.
\end{proposition}
\begin{proof}
  Let $\Theta$ and $\widetilde{\Theta}$ be two module traces on $\CM$. According to Theorem \ref{theorem:equivalences-for-C-trace} they correspond to $\Cat{C}$-balanced natural isomorphisms $\eta,\widetilde{\eta}: \Hom_{}(m,n) \rr \Hom_{}(n,m)^{*}$, respectively.
  Hence the vertical composition $\eta^{-1} \cdot \widetilde{\eta}: \Hom_{}(m,n) \rr \Hom_{}(m,n) $ of the natural isomorphisms
 is a $\Cat{C}$-balanced natural isomorphism.   According to Lemma \ref{lemma:basics-module-functor} there is  a unique $\Cat{C}$-module natural isomorphism 
  $Z: \id_{\Cat{M}} \rr \id_{\Cat{M}}$ such that 
  \begin{equation}
    \eta^{-1} \cdot \widetilde{\eta}(f)= Z(n) \circ f \quad \text{for all} \quad f \in \Hom_{}(m,n). 
  \end{equation}
  Theorem \ref{theorem:Dual-is-fusion} implies  that 
  there is a non-zero complex number $z$ such that $Z(f)=z \cdot f$ for all morphisms $f$ in $\Cat{M}$. Thus $\widetilde{\eta}(f)=z \cdot \eta(f)$ and so $\widetilde{\Theta}=z \cdot \Theta$.
\end{proof}

We refer to the transformation $eta \mapsto z \cdot \eta$ as rescaling of a module trace $\eta$ in the sequel.

\subsection{The Double Adjoints of Module Functors}

In this subsection we construct  module natural isomorphisms between   module functors of module categories with 
$\Cat{C}$-module traces and their double adjoint module functors. These isomorphisms are compatible with the composition of functors and if the module category is indecomposable
they define a  pivotal structure for the dual fusion category.
Recall that the left and right adjoint functors of a module functor $F: \CM \rr \CN$, $F^{l}$ and $F^{r}$, have a canonical structure of module functors. Note that in our convention the left adjoint functor $F^{l}$ is 
a right dual object to $F$ in the tensor category of functors and natural transformations. 
\begin{theorem}
  \label{theorem:module-functors-ambidextrous}
  Consider $\CM, \CN \in \ModCt$. For all  module functors $F:\CM \rr \CN$
  there is a canonical module  natural  isomorphism $a_{F}: F \rr F^{ll}$ 
  to the double left adjoint module functor of $F$.
  \begin{theoremlist}
  \item \label{item:natural-aF} The natural isomorphisms $a_{F}$ are natural with respect to module natural transformations, i.e.~for any module functor $G: \CM \rr \CN$  and 
    any module natural transformation $\rho: F \rr G$, the diagram 
    \begin{equation}
      \begin{xy}
        \xymatrix{ 
          F \ar[r]^{a_{F}} \ar[d]_{\rho}& F^{ll} \ar[d]^{\rho^{ll}} \\
          G   \ar[r]^{a_{G}}   &   G^{ll} 
        }
      \end{xy}
    \end{equation}
    commutes.
  \item  \label{item:comp-with-comp-aF}For all module functors $F:\CM \rr \CN$ and $K: \CN \rr \CE$,
    \begin{equation}
      \label{eq:compat-of-adj-comp}
      a_{KF}=a_{K}\circ a_{F}: K\circ F \rr (K \circ F)^{ll}.
    \end{equation}
  \end{theoremlist}
  In particular, these isomorphisms equip the dual category $\CMstar=\Funl{\Cat{C}}{\Cat{M},\Cat{M}}$ with a pivotal structure that is invariant under rescaling of the module trace of $\CM$. 
\end{theorem}
\begin{proof}
 According to Theorem \ref{theorem:equivalences-for-C-trace} we can identify the module traces with $\Cat{C}$-balanced natural isomorphisms $\eta^{M}:\Hom_{}(m,\tilde{m}) \rr \Hom_{}(\tilde{m},m)^{*}$ and 
  $\eta^{\Cat{N}}:\Hom_{}(n, \tilde{n}) \rr \Hom_{}(\tilde{n},n)^{*}$. 
  Consider the following sequence of natural $\Cat{C}$-balanced isomorphisms: 
  \begin{equation}
    \label{eq:left-right-identification}
    \begin{split}
      \Hom_{\Cat{N}}(n,F(m)) & \simeq  \Hom_{\Cat{M}}(F^{l}(n),m) \stackrel{\eta^{\Cat{M}}}{\simeq}  \Hom_{\Cat{M}}(m,F^{l}(n))^{*}\\
      & \simeq \Hom_{\Cat{N}}(F^{ll}(m),n)^{*}  \stackrel{(\eta^{\Cat{N}})^{-1}}{\simeq}\Hom_{\Cat{N}}(n, F^{ll}(m)).
    \end{split}
  \end{equation}
  According to  Lemma \ref{lemma:basics-module-functor}, the composition defines a  $\Cat{C}$-module natural isomorphism   $a_{F}: F \rr F^{ll}$.

  For the first part we have to show that the diagram
  \begin{equation}
    \begin{xy}
      \xymatrix{
        \Hom_{}(n, Fm) \ar[d]^{\simeq}              \ar[rr]^{\Hom(n,\rho m)}    \ar@/_4pc/[dddd]_{a_{F}}         &    &    \Hom_{}(n,Gm)  \ar[d]_{\simeq}  \ar@/^4pc/[dddd]^{ a_{G}}    \\
        \Hom_{}(F^{l}n,m)  \ar[d]^{\eta^{\Cat{M}}}       \ar[rr]^{\Hom(\rho^{l} n,m)}                                &    &   \Hom_{}(G^{l}n,m)       \ar[d]_{\eta^{\Cat{M}}}       \\
        \Hom_{}(m,F^{l}n)^{*}  \ar[d]^{\simeq}          \ar[rr]^{\Hom(m,\rho^{l}n)^{*}}                                &  &  \Hom_{}(m,G^{l}n)^{*}      \ar[d]_{\simeq}     \\
        \Hom_{}(F^{ll}m,n)^{*} \ar[d]^{(\eta^{\Cat{N}})^{-1}} \ar[rr]^{\Hom(\rho^{ll}m,n)^{*}}    &                  &  \Hom_{}(G^{ll}m,n)^{*}\ar[d]_{(\eta^{\Cat{N}})^{-1}} \\
        \Hom_{}(n,F^{ll}m)                                \ar[rr]^{\Hom(n,\rho^{ll}m)}       &             & \Hom_{}(n,G^{ll}m)
      }
    \end{xy}
  \end{equation}
  commutes.  All subdiagrams commute either by naturality of $\eta^{\Cat{M}}$ and $\eta^{\Cat{N}}$, by definition of the adjoint of $\rho$, or by definition of $a_{F}$ and $a_{G}$. Hence the 
  whole diagram commutes.

  For the second part we identify $(KF)^{l}=F^{l}K^{l}$. 
  It is enough to prove that the following diagram commutes:
  \begin{equation}
    \label{eq:adjoint-functor-comp}
    \begin{xy}
      \xymatrix{
        \Hom_{}(e,KFm)  \ar[d]^{\simeq} \ar[drr]^{a_{KF}} \ar@/_6pc/[ddddd]^{Ka_{F}}  &                  &                                \\
        \Hom_{}(F^{l}K^{l}e,m)        \ar[d]^{\eta^{\Cat{M}}}&                    &                       \Hom_{}(e,K^{ll}F^{ll}m)         \\
        \Hom_{} (m,F^{l}K^{l}e)^{*}\ar[d]^{\simeq}  &       \Hom_{}(K^{ll}F^{ll}m,e)^{*}\ar[ur]^{(\eta^{\Cat{E}})^{-1}}                     &                          \\
        \Hom_{}(F^{ll}m,K^{l}e)^{*}    \ar[d]^{(\eta^{\Cat{N}})^{-1}} \ar[ur]^{\simeq}     &                          &                        \\
        \Hom_{}(K^{l}e,F^{ll}m) \ar[d]^{\simeq}       &                           &                                \\
        \Hom_{}(e,KF^{ll}m)  \ar@/_4pc/[uuuurr]^{a_{K}F^{ll}}.   &                            &               
      }
    \end{xy}
  \end{equation}
  The upper triangle and the lower subdiagram commute due to the definition of $a_{KF}$ and $a_{K}$, respectively. It remains to show that the subdiagram on the left commutes. 
  It is easy to see that this subdiagram can be rewritten as 
  \begin{equation}
\label{eq:alternative-detail}
    \begin{xy}
      \xymatrix{
        \Hom_{}(e,KFm) \ar[r]^{\simeq}\ar[d]^{\Hom_{}(e,Ka_{F}m)}     &\Hom_{}(K^{l}e,Fm) \ar[d]^{\Hom_{}(K^{l}e,a_{F}m)} \\
        \Hom_{}(e,KF^{ll}m) \ar[r]^{\simeq}   & \Hom_{}(K^{l}e, F^{ll}m),
      }
    \end{xy}
  \end{equation}
 The commutativity of the diagram (\ref{eq:alternative-detail}) follows from the naturality of the adjunction and thus the second part is proven.
  From  part \refitem{item:natural-aF} and \refitem{item:comp-with-comp-aF} it is clear that the isomorphisms $a_{F}$ equip   $\CMstar=\Funl{\Cat{C}}{\Cat{M},\Cat{M}} \ni F $ with a pivotal structure. As the construction of $a_{F}$ involves 
  the map $\eta^{\Cat{M}}:\Hom_{\Cat{M}}(m,\tilde{m})\simeq \Hom_{\Cat{M}}(\tilde{m},m)^{*}$ composed with its inverse,  a constant scale factor cancels out.
\end{proof}

\begin{corollary}
  \label{corollary:piv-on-dual-module-trace}
  Let $\CM \in \ModCt$. Consider $\Cat{M}$ as a $\CMstar$-left module category and equip $\CMstar$ with  the induced pivotal structure from Theorem \ref{theorem:module-functors-ambidextrous}. 
  Then the $\Cat{C}$-module 
  trace on $\Cat{M}$ is also a $\CMstar$-module trace.
\end{corollary}
\begin{proof}
See Section \ref{sec:module-categories} for the structures of the category  $\CMstar=\Funl{\Cat{C}}{\Cat{M},\Cat{M}}$. The action of a $F \in \CMstar$ on $\CM$ is given by the application of the functor $F$. 
 By Theorem \ref{theorem:equivalences-for-C-trace} it is sufficient to show that the $\Cat{C}$-balanced natural isomorphism $\eta:\Hom_{}(m,\tilde{m}) \simeq \Hom_{}(\tilde{m},m)^{*}$ is also $\CMstar$-balanced. 
  The induced pivotal structure provides a natural isomorphism  $a_{F}^{r}: F^{r} \rr F^{l}$ for a functor $F \in \CMstar$.  
We have to show that the  diagram
  \begin{equation}
    \label{eq:also-CMStar-balanced}
    \begin{xy}
      \xymatrix{
        \Hom_{}(m,Fn) \ar[r]^{\eta} \ar[dd]^{\simeq} & \Hom_{}(Fn,m)^{*}\ar[d]^{\simeq}\\
        &  \Hom_{}(n,F^{r}m)^{*} \ar[d]^{(a_{F}^{r})^{-1}}\\
        \Hom_{}(F^{l}m,n) \ar[r]^{\eta} & \Hom_{}(n,F^{l}m)^{*}
      }
    \end{xy}
  \end{equation}
  commutes for all $m,n \in \Cat{M}$ and $F \in \CMstar$. The arrows pointing downwards are the $\Cat{C}$-
balancing natural isomorphism for $\Hom_{}(m,n)$ and $\Hom_{}(n,m)^{*}$, that are defined by the adjunction and the pivotal structure 
 according to equation (\ref{eq:balancing-varphi1}) and (\ref{eq:balancing-varphi2}), respectively. 
  The natural isomorphism $a_{F}$ is defined by equation (\ref{eq:left-right-identification}) in precisely such a way that the diagram commutes. Hence the statement follows.
 \end{proof}

\subsection{Conjugation of Pivotal Structures}
\label{sec:conj-pivot-struct}
When we restrict the considerations of the previous subsection  to the case of $\Cat{C}$ as a left module category over itself, Theorem \ref{theorem:module-functors-ambidextrous} provides a conjugation operation on the set of 
pivotal structures of a fusion category $\Cat{C}$. We show how this conjugation can be obtained alternatively from a canonical natural monoidal isomorphism $g: \id_{\Cat{C}} \rr (.)^{****}$ that exists for all 
fusion categories $\Cat{C}$.
\begin{theorem}
  \label{theorem:Dual-pivotal-structure}
Let $\Cat{C}$ be a fusion category with pivotal structure $a: \id_{\Cat{C}}\rr (.)^{**}$. 
\begin{theoremlist}
  \item \label{item:exist-conj}There exists a pivotal structure $\cc{a}:\id_{\Cat{C}}\rr (.)^{**}$ for $\Cat{C}$ with 
$(\cc{a}_{\leftidx{^{**}}{x}{}})^{-1}: x\rr \leftidx{^{**}}{x}{}$ defined by
\begin{equation}
\label{eq:dual-piv}
\ifx\du\undefined
  \newlength{\du}
\fi
\setlength{\du}{10\unitlength}

\;,
\end{equation}
for all $f \in \Hom_{}(c, d \otimes x)$ and $g \in \Hom_{}(d \otimes \leftidx{^{**}}{x}{},c)$.
\item The  dimension of an object $x$ with respect to the pivotal structure $\cc{a}$ is equal to the dimension of $\leftidx{^{*}}{x}{}$ with respect to   $a$. 
\item \label{item:condition-spherical}$\cc{a}=a$ if and only if $a$ is spherical.
\item \label{item:double-conj-equl}$\cc{\cc{a}}=a$.
\end{theoremlist}
\end{theorem}
\begin{proof}
It is well-known (see e.g. \cite{ENOfus}) that $\Funl{\Cat{C}}{\Cat{C},\Cat{C}}$ is canonically equivalent to $\rev{C}$ as a fusion category. $\rev{C}$ is the category $\Cat{C}$ with 
reversed tensor product. The module functors $\CC \rr \CC$ can be identified with 
the functors $(.) \otimes x$ of right tensoring with objects $x \in \Cat{C}$. The left adjoint functor to $(.) \otimes x$ is $(.) \otimes \leftidx{^{*}}{x}{}$.
To show \refitem{item:exist-conj}, we introduce the following graphical notation for the isomorphism $\eta^{\Cat{C}}:\Hom_{}(c,d) \rr \Hom_{}(d,c)^{*}$:
\begin{equation}
\ifx\du\undefined
  \newlength{\du}
\fi
\setlength{\du}{10\unitlength}

.
\end{equation}
Once the ellipse is replaced by a morphism  $ h \in \Hom_{}(d,c)$, the diagram represents the number $\eta^{\Cat{C}}(f)(h)$. The chain of isomorphisms (\ref{eq:left-right-identification}) reads now in 
graphical terms:
\begin{equation}
  \begin{split}
     &
\ifx\du\undefined
  \newlength{\du}
\fi
\setlength{\du}{10\unitlength}

  \end{split}
\end{equation}
Inserting once more the definition of $\eta^{\Cat{C}}$, we conclude that  equation (\ref{eq:left-right-identification}) yields equation (\ref{eq:dual-piv}). Hence Theorem \ref{theorem:module-functors-ambidextrous}
implies the first part.

The second statement follows by restricting the first statement to the case $d=\leftidx{^{*}}{x}{}$, $c=1$, $f=\coevp{x}$ and $g=\evp{\leftidx{^{*}}{x}{}}$. Recall that we defined the dimensions in a pivotal 
category as the left trace of the identity morphism.

Now consider the case  $a=\cc{a}$. The second part implies $\dim^{\Cat{C}}(c)= \dim^{\Cat{C}}(\leftidx{^{*}}{c}{})$ for all $c \in \Cat{C}$ and it follows that $a$ is spherical (see \cite{Mue1}).
Conversely, suppose that $a$ is spherical. Then
\begin{equation}
\ifx\du\undefined
  \newlength{\du}
\fi
\setlength{\du}{10\unitlength}

\;,
\end{equation}
where we used that $a$ is spherical in the last step. So $\cc{a}=a$ by equation (\ref{eq:dual-piv}). 
For part \refitem{item:double-conj-equl} we have to show that 
\begin{equation}
\ifx\du\undefined
  \newlength{\du}
\fi
\setlength{\du}{10\unitlength}

.
\end{equation}
With the symmetry of the left trace we calculate
\begin{equation}
\ifx\du\undefined
  \newlength{\du}
\fi
\setlength{\du}{10\unitlength}
%
\;,
\end{equation}
where in the last step we used equation (\ref{eq:dual-piv}) with the morphism $g$ in (\ref{eq:dual-piv}) set to  $\id_{\leftidx{^{**}}{x}{}}$. This proves the theorem.
\end{proof}
We call the pivotal structure $\cc{a}$ the conjugate pivotal structure of $a$. In the example  of $G$-graded vector spaces, see \ref{example:finite-G}, where a pivotal structure is a group homomorphism 
$\kappa:G \rr \C$, the conjugate pivotal structure is indeed given by the complex conjugate of $\kappa$.

It is instructive to consider the existence of conjugate pivotal structures also from another perspective. In \cite{ENOfus} it  is shown 
that for every fusion category there exists a monoidal natural isomorphism $g:\id \rr (.)^{****}$. We provide a simple description of such an isomorphism using dual $\Hom_{}$-spaces and show that the conjugate of a pivotal 
structure can be constructed with this isomorphism. 
We remark that in  \cite{HaggeHong} another graphical  proof of the existence of such a natural isomorphism $g$ is given with  a different approach to pivotal structures. 
\begin{proposition}
  \label{proposition:mon-iso-g}
Let $\Cat{C}$ be a fusion category. 
\begin{propositionlist}
  \item \label{item:natural-phi}The map 
    \begin{equation}
      \phi_{c}: \Hom_{}(c, 1) \rr \Hom_{}(1,c)^{*}, \quad \phi(f)(h)= h \circ f \in \C
    \end{equation}
for $c \in \Cat{C}$, $f \in \Hom_{}(1,c)$ and $h \in \Hom_{}(c,1)$ is a natural isomorphism.
\item \label{item:natural-g} The following chain of isomorphisms 
  \begin{equation}
\label{eq:definition-g}
    \begin{split}
        &\Hom_{}(x,\leftidx{^{**}}{c}{}) \simeq \Hom_{}(\leftidx{^{*}}{c}{} \otimes x,1) \stackrel{\phi}{\simeq} \Hom_{}(1, \leftidx{^{*}}{c}{} \otimes x)^{*} \simeq \Hom_{}(c,x)^{*} \\
& \simeq \Hom_{}(1,x \otimes c^{*})^{*} \stackrel{\phi^{-1}}{\simeq} \Hom_{}(x \otimes c^{*},1) \simeq \Hom_{}(x, c^{**})
    \end{split}
    \end{equation}
is natural in $c,x \in \Cat{C}$ and defines a monoidal natural isomorphism $g_{c}: \leftidx{^{**}}{c}{} \rr c^{**}$.
\item  \label{item:graphical-g} $g_{c}:\leftidx{^{**}}{c}{} \rr c^{**}$ is  defined  uniquely by the requirement that for all $f \in \Hom_{}(x, \leftidx{^{**}}{c}{})$ and $h\in \Hom_{}(c,x)$:
  \begin{equation}
    \label{eq:def-g-graph}
\ifx\du\undefined
  \newlength{\du}
\fi
\setlength{\du}{10\unitlength}

\;.
  \end{equation}
\end{propositionlist}
\end{proposition}
\begin{proof}
  The naturality of $\phi$ in part \refitem{item:natural-phi} is clear.   $\phi$ is an isomorphism due to the semisimplicity of $\Cat{C}$. For part \refitem{item:natural-g}, the naturality  
 of the isomorphisms in 
$x$ and $c$ is a consequence of part \refitem{item:natural-phi} and the naturality of the duality.
Hence the isomorphism $g_{c}$ is well-defined  by the Yoneda lemma.  
We introduce the graphical notation 
\begin{equation}
  \label{eq:graph-not-phi}
\ifx\du\undefined
  \newlength{\du}
\fi
\setlength{\du}{10\unitlength}

\end{equation}
for $\phi_{c}(f) \in \Hom_{}(1,c)^{*}$. If the  unlabelled ellipse  is replaced by an morphism  $ h \in \Hom_{}(1,c)$, this expression represents the number $\phi_{c}(f)(h)$.
Now the chain of isomorphisms (\ref{eq:definition-g}) reads in graphical terms
\begin{equation}
  \begin{split}
 &
\ifx\du\undefined
  \newlength{\du}
\fi
\setlength{\du}{10\unitlength}

,
  \end{split}
\end{equation}
where   $\tilde{f}$ is defined by
\begin{equation}
\ifx\du\undefined
  \newlength{\du}
\fi
\setlength{\du}{10\unitlength}
%
.
\end{equation}
Applying the rigidity of $\Cat{C}$ it follows that 
\begin{equation}
\label{eq:almost-there-graph-phi}
\ifx\du\undefined
  \newlength{\du}
\fi
\setlength{\du}{10\unitlength}
%
.
\end{equation}
Applying  once more the rigidity of $\Cat{C}$, equation (\ref{eq:almost-there-graph-phi}) implies expression (\ref{eq:definition-g}).
For the compatibility of $g$ with the monoidal structure we calculate
\begin{equation}
  \label{eq:1}
\ifx\du\undefined
  \newlength{\du}
\fi
\setlength{\du}{10\unitlength}

\; ,
\end{equation}
where we first used the graphical expression (\ref{eq:def-g-graph}) for $g_{d}$, then for $g_{c}$ and finally for $g_{d \otimes c}$.
Since this equality holds for all morphisms $h \in \Hom_{}(c \otimes d, x)$ and $f \in \Hom(x,  \leftidx{^{**}}{d}{} \otimes  \leftidx{^{**}}{c}{})$, we conclude that $g_{d \otimes c}= g_{d} \otimes g_{c}$ from the uniqueness statement in part \refitem{item:graphical-g}.
\end{proof}
\begin{remark} \footnote{We are grateful to the referee for bringing \cite{EtinAnalogue} to our attention.}
In \cite[Theorem 7.3]{EtinAnalogue}, a canonical monoidal isomorphism $\delta: (.)^{**} \simeq  {}^{**}(.)$   is defined by $\tr(\phi \circ \delta_{c}^{-1})=\tr(\phi^{*})$ 
for all isomorphisms $\phi: c^{**} \simeq c$ for a simple object $c$ in a fusion category.
We show that $\delta_{c}^{-1}$ coincides with $g_{c}$ as defined by Proposition \ref{proposition:mon-iso-g}.
Let  $\phi: c^{**} \simeq c$ be an isomorphism for a simple object c. Then 
\begin{equation}
  \label{eq:2}
  \tr(\phi \circ g_{c})=
\ifx\du\undefined
  \newlength{\du}
\fi
\setlength{\du}{10\unitlength}

=\tr(\phi^{*}) 
\end{equation}
shows that our definition of the natural isomorphism  $g$ coincides with the definition in \cite{EtinAnalogue} and hence also with the definition of $g$ in\cite{ENOfus}. The advantage of our definition is that it is defined directly for all objects and not just for simple objects.
\end{remark}

The following proposition  clarifies the relation between $g$ and the conjugate of a pivotal structure.
\begin{proposition}
  Let $\Cat{C}$ be a fusion category with pivotal structure $a: \id \rr (.)^{**}$. 
  \begin{propositionlist}
    \item a and its conjugate $\cc{a}$ combine to $g$, i.e.~$ \cc{a}_{c} \cdot a_{\leftidx{^{**}}{c}{}}= a_{c} \cdot \cc{a}_{\leftidx{^{**}}{c}{}} =g_{c}: \leftidx{^{**}}{c}{} \rr c^{**}$. 
\item $a$ is spherical if and only if $a_{c} \cdot a_{\leftidx{^{**}}{c}{}}=g_{c}$.
  \end{propositionlist}
\end{proposition}
\begin{proof}
  For all $f:c \rr {\leftidx{^{**}}{c}{}}$, 
  \begin{equation}
\ifx\du\undefined
  \newlength{\du}
\fi
\setlength{\du}{10\unitlength}

\;,
  \end{equation}
by equation (\ref{eq:dual-piv}). This implies $ a_{c} \cdot \cc{a}_{\leftidx{^{**}}{c}{}} =g_{c}$ with condition (\ref{eq:def-g-graph}).  The other equation follows directly from the naturality of $a$.
For the second part note that the first part implies $ \cc{a}_{c}=g_{c} \cdot  a_{\leftidx{^{**}}{c}{}}^{-1}$. Now the statement follows directly from Theorem \ref{theorem:Dual-pivotal-structure}, \refitem{item:condition-spherical}. 
The statement can also be derived directly from the graphical  expression (\ref{eq:def-g-graph}). 
\end{proof}

\section{The Existence Problem as an Eigenvalue Equation}
\label{sec:existence-problem-as}
The aim of this section is to formulate the existence of a module trace as an eigenvalue problem. In particular this allows  to conclude from known results in the literature, that all module categories over pseudo-unitary fusion categories equipped with 
the canonical spherical structure admit a module trace. 
\subsection{The Dimension Matrix of a Module Category}
We show how a trace on a semisimple category is characterised by the dimensions of simple objects using the trace in $\Vect$. For a module trace on a module category over $\Cat{C}$ we derive an analogous formula with 
the trace in $\Vect$ replaced by 
the left trace in $\Cat{C}$.
As a consequence we obtain that the existence of  a module trace on $\CM$ implies $\dim^{\Cat{C}}(\icm{m,m}) > 0$
for all simple $m \in \Cat{M}$. In the last part we apply the considerations to spherical fusion categories and show that a pivotal structure for $\Cat{C}$ is spherical if and only if there is a module category $\CM$
over $\Cat{C}$ with a module trace such that all dimensions in $\Cat{M}$ are real.

Consider general traces on a semisimple category $\Cat{M}$ with a finite set of representatives $m_{i}$, $i \in I$ for the isomorphism classes of simple objects. 
 The following lemma is well-known, see e.g. \cite[Lemma \Rmnum{2}.4.2.3]{Tur}.
\begin{lemma}
  \label{lemma:general-trace-ses}
  A collection of linear maps $\Theta_{m}: \End_{\Cat{M}}(m) \rr \C$ that satisfies the symmetry property of Definition \ref{definition:Module-trace} \refitem{item:Theta-symmetric} is non-degenerate 
  and hence a trace on $\Cat{M}$ if and only if $\Theta(\id_{m_{i}}) \neq 0$ for all $i \in I$.
\end{lemma}
\begin{proposition}
 For every   trace $\Theta$  on $\Cat{M}$,  $ (\dim^{\Theta}(m_{i}))_{i \in I}$ is an $|I|$-tuple of non-zero numbers. Conversely, given such a tuple $d_{i}\in \C^{\times}$, $i \in I$, 
\begin{equation}
    \label{eq:alt-formula-theta}
    \Theta_{m}(f)= \sum_{i \in I} \tr(\Hom_{}(m_{i},f))d_{i}, 
  \end{equation}
  for $f \in \Hom_{}(m,m)$ defines a trace on $\Cat{M}$. Here $\tr(\Hom_{}(m_{i},f))$ denotes 
  the usual trace on $\Vect$ of the linear map $\Hom(m_{i},f): \Hom(m_{i},m) \rr \Hom(m_{i},m)$.

These two maps  yield a bijection between the set of  traces  on $\Cat{M}$ and the set of  $|I|$-tuples of non-zero numbers. 
\end{proposition}
\begin{proof}
  Suppose that $\Cat{M}$ is equipped with  a trace $\Theta$. Then   $d_{i}=\dim^{\Theta}(m_{i}) \neq 0$
  due to Lemma \ref{lemma:general-trace-ses}. We have to show that for all $f \in \End(m)$ 
  formula (\ref{eq:alt-formula-theta}) holds. 
  The semisimplicity of $\Cat{M}$ ensures that the functor 
  \begin{equation}
    \label{eq:semisimple-functor}
    \Cat{M} \ni m \mapsto  \oplus_{i}\Hom_{\Cat{M}}(m_{i},m) \tensor[\C] m_{i} 
  \end{equation}
  is naturally isomorphic to the identity functor on $\Cat{M}$.
  This implies  
  \begin{equation}
    \begin{split}
      \Theta_{m}(f)&= \Theta_{\oplus_{i}\Hom_{}(m_{i},m) \tensor[\C] m_{i}}(\oplus_{i} \Hom(m_{i},f) \tensor[\C] m_{i})\\
      &= \sum_{i \in I} \Theta_{\Hom_{}(m_{i},m) \tensor[\C] m_{i}}(\Hom(m_{i},f) \tensor[\C] m_{i}) \\
      &=\sum_{i \in I} \tr(\Hom_{}(m_{i},f))d_{i},
    \end{split}
  \end{equation}
  where we used Lemma \ref{lemma:Vect-mdoule-trace} in the last step. 

  For the converse we have to show that  given a set of non-zero $d_{i} \in \C$ for $i\in I$, formula (\ref{eq:alt-formula-theta}) defines a trace on $\Cat{M}$.  The symmetry follows directly from the cyclic property of $\tr$.  
  The non-degeneracy follows from Lemma \ref{lemma:general-trace-ses}.
\end{proof}

Now we discuss  $\Cat{C}$-module traces.
First we need a technical result. Choose representatives $c_{u}$, $u \in U$ for the isomorphism classes of simple objects of $\Cat{C}$.
  See e.g. \cite{Runkel} for a review of the definition of the Deligne product $\boxtimes$ of additive categories.
\begin{lemma}
  \label{lemma:special-isomorphic-functors}
  The following functors $ \Cat{M} \rr \Cat{C} \boxtimes \Cat{M}$ are naturally isomorphic. 
  \begin{equation}
    \label{eq:M-box-D-to-D}
    \begin{split}
      m&   \mapsto \oplus_{u \in U} c_{u} \boxtimes  \leftidx{^{*}}{c_{u}}{} \act m, \quad \text{and}\\
      m & \mapsto \oplus_{i \in I} \icm{m_{i},m} \boxtimes m_{i}.
    \end{split}
  \end{equation}
\end{lemma}

\begin{proof}
  The objects $ \oplus_{u \in U} c_{u} \boxtimes  \leftidx{^{*}}{c_{u}}{} \in \Cat{C} \boxtimes \Cat{C}$ and 
  $\oplus_{i \in I} m_{i} \boxtimes m_{i} \in \op{M} \boxtimes \Cat{M}$ are independent of the choice of  representatives of simple objects in the sense that  the objects obtained from any two
  choices of simple objects are canonically isomorphic, see \cite[Sec. 2.4]{BakKir}. This shows that the two maps yield well-defined functors. 
  Now let $c \boxtimes \tilde{m} \in \Cat{C}\boxtimes \Cat{M}$. Using the semisimplicity of $\Cat{C}$ and $\Cat{M}$ we obtain the following chain of natural isomorphisms:
  \begin{equation}
    \label{eq:rhs-hom-phi1-phi2}
    \begin{split}
      \Hom_{\Cat{C}\boxtimes\Cat{M}}(c\boxtimes \tilde{m}, \oplus_{u}  c_{u}\boxtimes \leftidx{^{*}}{c_{u}}{}\act {m}) &  \simeq \oplus_{u}  \Hom_{\Cat{C}}(c,c_{u}) \otimes  \Hom_{\Cat{M}}(\tilde{m},  \leftidx{^{*}}{c_{u}}{} \act m) \\
      &  \simeq  \oplus_{u} \Hom_{\Cat{C}}(c,c_{u}) \otimes \Hom_{\Cat{C}}(c_{u},\icm{\tilde{m},m})\\
      & \simeq \Hom_{\Cat{C}}(c, \icm{\tilde{m},m}) \simeq \Hom_{}(c \act \tilde{m},m)\\
      &\simeq \Hom_{\Cat{M}}( \tilde{m},  \leftidx{^{*}}{c}{} \act m)\\
      &\simeq  \oplus_{i} \Hom_{\Cat{M}}(m_{i},  \leftidx{^{*}}{c}{} \act m) \otimes \Hom_{\Cat{M}}(\tilde{m},m_{i})\\
      &\simeq  \oplus_{i} \Hom_{\Cat{C}}(c, \icm{m_{i},m}) \otimes  \Hom_{\Cat{M}}(\tilde{m},m_{i}) \\
      &\simeq \Hom_{\Cat{C}\boxtimes \Cat{M}}(c \boxtimes \tilde{m}, \oplus_{i} \icm{m_{i},m} \boxtimes m_{i}).
    \end{split}
  \end{equation}
  Now apply the Yoneda lemma to obtain a natural isomorphism between the two functors. 
\end{proof}
The following result provides an  alternative characterisation of module traces.
Recall from \cite{ENOfus} that for a pivotal fusion category $\dim(\Cat{C})=\sum_{u \in U} |\dim^{\Cat{C}}(c_{u})|^{2} \neq 0$. 
\begin{proposition}
  \label{proposition:eigenvalue-problem-module-trace-ses}
  Let $\CM$ be a $\Cat{C}$-module category. 
  If $\Theta$ is a $\Cat{C}$-module trace on $\Cat{M}$, the dimension vector $d_{i}=\dim^{\Theta}(m_{i})$ for $i \in I$ consists of non-zero numbers $d_{i}$ and
  is a (right) eigenvector of the matrix $(Q)_{ij}=\dim^{\Cat{C}}(\icm{m_{j},m_{i}})$ with eigenvalue $\dim(\Cat{C})$. 
  If a tuple of non-zero numbers $d_{i}$ for $i \in I$ is an eigenvector of $(Q)_{ij}$ with eigenvalue $\dim(\Cat{C})$, 
  then the collection of linear maps
  \begin{equation}
    \label{eq:formula-C-trace}
    \Theta_{m}(f)= \frac{1}{\dim(\Cat{C})}  \sum_{i \in I} \tr^L(\icm{m_{i},f}) d_{i},
  \end{equation}
  for $f \in \End(m)$ and $m \in \Cat{M}$ defines a $\Cat{C}$-module trace on $\Cat{M}$. 
  These two maps are mutually inverse.
\end{proposition}

\begin{proof}
  Let $\Theta$ be $\Cat{C}$-module trace on $\Cat{M}$. Lemma \ref{lemma:special-isomorphic-functors}  implies that  
  the object $\oplus_{u} (c_{u} \otimes  \leftidx{^{*}}{c_{u}}{}) \act m$  is isomorphic to $\oplus_{i} \icm{m_{i},m}\act m_{i}$ in $\Cat{M}$. Hence,
  \begin{equation}
    \label{eq:dimension-eigenvector-Q}
    \begin{split}
      \dim(\Cat{C})\cdot d_{k}&=\dim^{\Cat{C}}(\oplus_{u} (c_{u} \otimes  \leftidx{^{*}}{c_{u}}{})) \cdot \dim^{\Theta}(m_{k})\\
      & =\dim^{\Theta}(\oplus_{u} (c_{u} \otimes  \leftidx{^{*}}{c_{u}}{}) \act m_{k})=\dim^{\Theta} (\oplus_{i} \icm{m_{i},m_{k}}\act m_{i}) \\
      & =\sum_{i \in I} \dim^{\Cat{C}}(\icm{m_{i},m_{k}}) d_{i}.
    \end{split}
  \end{equation}
  In the sequel  we will refer to the matrix $Q=(Q_{ij})$ as the dimension matrix  and to the vector $d=(d_{i})$ as the dimension vector. Equation (\ref{eq:dimension-eigenvector-Q}) shows that the dimension vector is a right eigenvalue 
of the dimension matrix with eigenvalue $\dim(\Cat{C})$. 
  As another consequence of Lemma \ref{lemma:special-isomorphic-functors} we obtain the identity
  \begin{equation}
    \Theta( \oplus_{u} (c_{u} \otimes  \leftidx{^{*}}{c_{u}}{}) \act f)= \Theta ( \oplus_{i} \icm{m_{i},f} \act m_{i}),
  \end{equation}
 for all $f \in \End(m)$. This implies  formula (\ref{eq:formula-C-trace})  
  with $d_{i}=\dim^{\Theta}(m_{i})$.

  Now suppose we are given an eigenvector $d$ of the dimension matrix with eigenvalue $\dim(\Cat{C})$ whose components do not vanish. 
  Then define a linear map $\Theta_{m}: \End(m) \rr \C$ by the  formula (\ref{eq:formula-C-trace}). The symmetry of $\Theta$ follows from the cyclic property of the left trace $\tr^{L}$ of $\Cat{C}$. 
Since  $\Theta_{m_{i}}(\id_{m_{i}})=\frac{1}{\dim(\Cat{C})}  \sum_{j}Q_{ij} d_{j}=d_{i}\neq 0$,  we conclude 
  with Lemma \ref{lemma:general-trace-ses} that $\Theta$ is a trace on $\Cat{M}$.  For the $\Cat{C}$-compatibility we have to show that for all $ f \in \End(c \act m)$, 
  \begin{equation}
    \label{eq:C-comp-of-formula}
    \sum_{i \in I} \tr^L(\icm{m_{i},f}) d_{i}= \sum_{i \in I} \tr^{L}(\icm{m_{i},\tr^{\Cat{C}}(f)})d_{i}. 
  \end{equation}
  Since $\icm{m_{i},.}: \CM \rr \CC$ is a module functor, Lemma \ref{lemma:partial-trace-module-funct} implies that $\tr^{\Cat{C}}(\icm{m_{i},f})=\icm{m_{i},\tr^{\Cat{C}}(f)}$. 
  Now the statement follows from $\tr^L(\icm{m_{i},f})=\tr^{L} \circ \tr^{\Cat{C}}(\icm{m_{i},f})$. 
\end{proof}
Note that formula (\ref{eq:formula-C-trace}) is a generalisation of formula (\ref{eq:alt-formula-theta}). 
\begin{remark}
  The proof of Proposition \ref{proposition:eigenvalue-problem-module-trace-ses} shows that for any set of numbers $d_{i}$, $i \in I$, formula (\ref{eq:formula-C-trace}) defines a 
  collection of linear maps that satisfy the symmetry and $\Cat{C}$-compatibility condition of Definition \ref{definition:Module-trace}. The non-degeneracy condition is fulfilled 
  if and only if $\sum_{j}Q_{ij}d_{j} \neq 0$ for all $i \in I$.
\end{remark}

Next we discuss some properties of the dimension matrix for a module category $\CM$ that not necessarily possesses a module trace.
Let
\begin{equation}
  \label{eq:notation}
  M_{u,i}^{j}=\dim_{\C}(\Hom_{\Cat{M}}(c_{u}\act m_{i},m_{j}))
\end{equation}
be the multiplicity matrix of the action of $c_{u}\in \Cat{C}$ on $\Cat{M}$.
\begin{proposition}\label{proposition:Q} Let $\CM$ be a $\Cat{C}$-module category.
  The dimension matrix $Q$ satisfies:
  \begin{propositionlist}
  \item \label{item:Qij} $Q_{ij}=\sum_{u \in U}\dim^{\Cat{C}}(c_{u}) M_{u,j}^{i}$.
  \item \label{item:Q-square} $Q^{2}= \dim(\Cat{C})\cdot  Q$.
  \item  $Q$ is hermitian.
  \end{propositionlist}
\end{proposition}

\begin{proof}
  The multiplicity of each object $c_{u}$ in $\icm{m_{j},m_{i}}$ is 
  \begin{equation}
    \dim_{\C}(\Hom_{}(c_{u}, \icm{m_{j},m_{i}}))=\dim_{\C}(\Hom_{}(c_{u}\act m_{j},m_{i}))=M_{u,j}^{i}.
  \end{equation}
 This shows part \refitem{item:Qij}.
 For the second claim we first compute
  \begin{equation}
    \label{eq:calcul-phi1-2-in-innerhom}
    \begin{split}
      \oplus_{j \in I}\icm{m_{j},m_{i}}\otimes \icm{m_{k},m_{j}}&=\oplus_{j \in j} \icm{ m_{k}, \icm{m_{j},m_{i}}\act m_{j}}\\
      &\simeq \icm{ m_{k}, \oplus_{u \in U}(c_{u} \otimes  \leftidx{^{*}}{c_{u}}{})\act m_{j}} \\
      &\simeq \oplus_{u \in U}(c_{u} \otimes  \leftidx{^{*}}{c_{u}}{})\otimes \icm{m_{k},m_{j}},
    \end{split}
  \end{equation}
  where we used Lemma \ref{lemma:special-isomorphic-functors} in the second step. Now the statement follows after applying $\dim^{\Theta}$ to both sides of this equation.
  For the third statement we show that the objects $\icm{m_{i},m_{j}}$ and $\icmstar{m_{j},m_{i}}$ are isomorphic in $\Cat{C}$. 
  We compute the multiplicity of a $c \in \Cat{C}$ in $\icmstar{m_{j},m_{i}}$ by using that $\Hom_{}(m,n) \simeq \Hom_{}(n,m)$ as vector spaces. The following isomorphisms are  isomorphisms of vector spaces:
  \begin{equation}
    \begin{split}
      \Hom_{}(c,  \icmstar{m_{j},m_{i}}) &\simeq\Hom_{}(\icm{m_{j},m_{i}}, \leftidx{^{*}}{c}{}) \simeq \Hom_{}( \leftidx{^{*}}{c}{}, \icm{m_{j},m_{i}})\\
      &\simeq \Hom_{}( \leftidx{^{*}}{c}{} \act m_{j},m_{i}) \simeq \Hom_{}(m_{i}, \leftidx{^{*}}{c}{} \act m_{j}) \\
      &\simeq \Hom_{}(c \act m_{i},m_{j})= \Hom_{}(c, \icm{m_{i},m_{j}}).
    \end{split}
  \end{equation}
As  the multiplicities of all simple objects agree, we conclude that there exists an isomorphism $\icm{m_{i},m_{j}} \rr \icmstar{m_{j},m_{i}}$ in $\Cat{C}$ .
 With  $\dim^{\Cat{C}}(c^{*})=\cc{\dim^{\Cat{C}}(c)}$ for all 
  objects $c \in \Cat{C}$ from \cite[Proposition 2.9]{ENOfus},  it follows that
  \begin{equation}
 \dim^{\Theta}(\icm{m_{i},m_{j}})= \dim^{\Theta}( \icmstar{m_{j},m_{i}}) =\cc{\dim^{\Theta}(\icm{m_{j},m_{i}})}.
     \end{equation}
\end{proof}

\begin{proposition}
  \label{proposition:Q-rank-1}
  A module category $\CM$ has a module trace if and only if the dimension matrix $Q$ is of rank $1$ with only non-zero entries. 
  In particular is $\dim^{\Cat{C}}(\icm{m,m}) >0$ for all simple objects $m \in \Cat{M}$.
\end{proposition}
\begin{proof}
  It follows  directly from  Proposition \ref{proposition:Q}, that the only possible (right and left) eigenvalues of $Q$ are $\dim(\Cat{C})$ and $0$. 
 Suppose $\Cat{M}$ has a module trace and $d$ is the corresponding  eigenvector of $Q$ with all entries non-zero.
  Let $\tilde{d}$ be an  eigenvector of $Q$  with eigenvalue  $\dim(\Cat{C})$. 
  There always exists a linear combination $d + \lambda \tilde{d}$ with all entries non-zero. Hence $\tilde{d}$ must be proportional to $d$. This shows that $Q$ has rank $1$.

 Now suppose $Q_{ij}=d_{j} \cc{d_{i}} $ with non-zero numbers $d_{i}$. Then $\sum_{i}\cc{d_{i}}d_{i}=\dim(\Cat{C})$ by Proposition \ref{proposition:Q} \refitem{item:Q-square}. Hence $d_{i}$ yields a module trace.
This proofs also the last statement since $\dim^{\Cat{C}}(\icm{m_{i},m_{i}})=d_{i}\cc{d_{i}}$.
\end{proof}
As an example we discuss pseudo-unitary fusion categories. Recall from \cite{ENOfus} the definition of the Frobenius-Perron dimensions of objects in a fusion category. 
A  pseudo-unitary fusion category possesses a canonical spherical structure such that the dimension of all objects are equal to the Frobenius-Perron dimensions. 
In  \cite{ENOfus} and \cite{ENOfuhom} Frobenius-Perron dimensions of simple objects in module categories are defined and studied. 
The following statement shows that for module categories over pseudo-unitary fusion categories our definition of 
module trace reduces to the Frobenius-Perron dimension of  \cite[Rem. 2.3]{ENOfuhom}.  
\begin{proposition}
  \label{proposition:pseudo-unitary}
  Let $\Cat{C}$ be a pseudo-unitary fusion category and let $\CM$ be a module category over $\Cat{C}$. The Frobenius-Perron dimensions of simple objects in $\Cat{M}$ provide a module 
trace for $\Cat{M}$ and thus is the canonical spherical structure of $\Cat{C}$ flexible. 
\end{proposition}
\begin{proof}
According to \cite{ENOfuhom},
  there exists a  Frobenius-Perron eigenvector  $(d_{i})_{i \in I}$ of $\Cat{M}$,  that is defined by  $d_{j}>0$ for all $j \in I$  and: 
  \begin{equation}
    \label{eq:explicit-FP-inM}
    \sum_{u \in U} M_{u,i}^{j}d_{j}=\dim^{\Cat{C}}(c_{u})d_{i}.
  \end{equation}
  If we multiply this equation with $\dim^{\Cat{C}}(c_{u})$, sum over $u \in U$ and use that the pivotal structure is spherical,
  we see that $(d_{i})$ is an eigenvector of $Q_{ij}$ with eigenvalue $\dim(\Cat{C})$ and hence defines a module trace according to Proposition \ref{proposition:eigenvalue-problem-module-trace-ses}.
\end{proof}

\subsection{Module Traces on Module Categories over Spherical Fusion Categories}
Next we discuss the relation of module traces and  spherical structure.

\begin{proposition}
  \label{proposition:can-assume-dimensions-real}
  Let $\Cat{C}$ be spherical, $\Cat{M}$ a left $\Cat{C}$-module category with module trace $\Theta$. 
  There exists a $ z \in \C$ such that the  dimensions of objects in $\Cat{M}$ with respect to the module trace $z \cdot \Theta$ are real.
\end{proposition}
\begin{proof}
  If $\Cat{C}$ is spherical all dimensions of $\Cat{C}$ are real. 
  Hence $Q$ is a real symmetric matrix which can be diagonalised by a real matrix. 
  If follows that  the entries of all eigenvectors of  $Q$ are real.
\end{proof}

The next result provides a criterion to determine whether  a given pivotal structure is spherical.
\begin{proposition}
  Let $\CM$ be a module category with module trace $\Theta$. 
  \begin{propositionlist}
  \item \label{item:dim-left}The dimension vector $d_{i}=\dim^{\Theta}(m_{i})$ is a left  
    eigenvector of the dimension matrix with eigenvalue $C=\sum_{u \in K} \dim (c_{u})^{2}$, i.e.~$\sum_{j} d_{j} Q_{ji}= C \cdot  d_{i}$.
  \item  \label{item:cirterion-C}The number $C=\sum_{u \in K} \dim (c_{u})^{2}$  is  equal to $\dim(\Cat{C})$ if and only if the pivotal structure is spherical, and it is equal to $0$ otherwise.
  \item \label{item:cirterion} A pivotal structure for $\Cat{C}$ is spherical  if and only if there exists a module category $\CM$ with module trace such that all dimensions of objects in $\Cat{M}$ are real.
  \item Let $\Cat{C}$ be spherical and assume that $\CM$ has a module trace. Then the induced pivotal structure for the dual category $\CMstar$ from Theorem \ref{theorem:module-functors-ambidextrous} is spherical.
  \end{propositionlist}
\end{proposition}
\begin{proof}
  The  $\Cat{C}$-compatibility of $\Theta$
  implies
  \begin{equation}
    \label{eq:action-dimensions-M}
    \sum_{i} M_{u,j}^{i}d_{i}=\dim^{\Theta}(c_{u}\act m_{j})= \dim^{\Cat{C}}(c_{u}) \cdot d_{j}.
  \end{equation}
  Multiplying this equation with $\dim^{\Cat{C}}(c_{u})$ and summing over $u \in U$ yields:
  \begin{equation}
    \label{eq:eigenvalue-Q-m}
    \begin{split}
      C \cdot d_{j} &= \sum_{i \in I, u \in U} \dim^{\Cat{C}}(c_{u}) M_{u,j}^{i}d_{i}\\
      &= \sum_{i \in I} \dim^{\Cat{C}}(\icm{m_{j},m_{i}}) d_{i},
    \end{split}
  \end{equation}
  where we used Proposition \ref{proposition:Q},~\refitem{item:Qij}. This proves the first statement.  For the module category $\CC$, 
  equation (\ref{eq:eigenvalue-Q-m}) implies  $ C\cdot \dim^{\Cat{C}}(c_{j})=C \cdot \dim^{\Cat{C}}(c_{j}^{*})$.  It is shown in \cite{Mue1} that  $\Cat{C}$ is spherical if and only if  
  $\dim^{\Cat{C}}(c_{u}^{*})=\dim^{\Cat{C}}(c_{u})$. Hence the second statement follows. To prove part \refitem{item:cirterion}, let $\CM$ be a module category with module trace $\Theta$ and $d_{i}= \dim^{\Theta}(m_{i}) \in \R$ for 
  all simple $m_{i} \in \Cat{M}$. According to Proposition \ref{proposition:Q-rank-1} we can assume, using the freedom to rescale $\Theta$, 
that $\sum_{i}d_{i}^{2}= \dim(\Cat{C})$ and therefore $Q_{ij}=d_{i} d_{j}$. From part \refitem{item:dim-left} it follows that 
  $C= \dim(\Cat{C})$ and hence part \refitem{item:cirterion-C}  implies that the pivotal structure is spherical. The converse is clear since the  module category $\CC$ has real dimensions if $a$ is spherical.
  The last statement is a consequence of part \refitem{item:cirterion} together with Proposition \ref{proposition:can-assume-dimensions-real} and Corollary \ref{corollary:piv-on-dual-module-trace}.
\end{proof}
\begin{remark}
  It is shown in \cite[Theorem 5.16]{Mue1} by different methods that an indecomposable 
module category $\CM$ over a spherical category $\Cat{C}$ provides a spherical structure for the category $\Funl{\Cat{C}}{\Cat{M},\Cat{M}}$.
  The relation to our construction remains to be investigated.
\end{remark}

\section{Frobenius Algebras}
\label{subsec:Frobenius-Algebras}
In this section  we show that  module traces are directly related to Frobenius algebras. This is done by exploring the graphical calculus for module categories with module traces and constructing 
a natural isomorphism  $\beta$ that is the reflected analogue of the $\alpha$ in Subsection \ref{subsec:Inner-hom}. This operation equips the inner hom objects with the structure of a Frobenius algebra. 
We also prove the converse, namely that the module category formed by the  modules over a special haploid symmetric Frobenius algebra  has a module trace.

To emphasise the role of the $\Cat{C}$-compatibility of a module trace we first discuss traces on a module category $\CM$.
We saw in Section \ref{subsection:Characterization-module-trace} that a module category with a trace that is not necessarily $\Cat{C}$-compatible 
equips $\CM$ with a  natural isomorphism $\eta^{\Cat{M}}: \Hom_{\Cat{M}}(m,n) \rr \Hom_{\Cat{M}}(n,m)^{*}$. 
Recall that the pivotal structure of $\Cat{C}$ also yields a trace and  a natural isomorphism 
$\eta^{\Cat{C}}: \Hom_{\Cat{C}}(c,d) \rr \Hom_{\Cat{C}}(d,c)^{*}$, see equation (\ref{eq:eta-C}).

\begin{proposition}
  \label{proposition:def-beta}
  Let $\CM$ be a $\Cat{C}$-module category equipped with a trace $\Theta$. Then  there exists a  natural isomorphism 
  \begin{equation}
    \beta:  \Hom_{}(n,c \act m) \rr   \Hom_{}(\icm{m,n},c),
  \end{equation}
which is specified uniquely by the requirement
  \begin{equation}
    \label{eq:specify-beta}
    \tr^{L}(\beta(f) \circ \alpha(g))= \Theta_{c \act m}(f \circ g), 
  \end{equation}
 for all $g \in \Hom_{}(c \act m,n)$ and with $f \in \Hom_{}(n,c \act m)$.
\end{proposition}
\begin{proof}
  Condition (\ref{eq:specify-beta}) is equivalent to  defining $\beta$ as the following composition  of natural isomorphisms:
  \begin{equation}
    \label{eq:inner-hom-both}
    \begin{split}
      & \Hom_{}( n,c \act  m) \stackrel{\eta^{\Cat{M}}}{\simeq} \Hom_{}(c \act m,  n)^{*} \\
      &\simeq \Hom_{}(c, \icm{m,n})^{*} \stackrel{(\eta^{\Cat{C}})^{-1}}{\simeq} \Hom_{}(\icm{m,n},c).
    \end{split}
  \end{equation}
  This follows directly from the identity $\eta^{\Cat{M}} (a)(b)=\Theta_{n}(a \circ b)$ for $a \in \Hom_{}(m,n)$ and $b \in \Hom_{}(n,m)$.
\end{proof}
The graphical representation of  $\beta$ is 
\begin{equation}
  \beta \corres 
\ifx\du\undefined
  \newlength{\du}
\fi
\setlength{\du}{10\unitlength}

\;,
\end{equation}
i.e.~$\beta$ allows one to flip strings representing objects in the module category upwards.
Equation (\ref{eq:specify-beta}) reads in graphical terms, where we omit the pivotal structure for better readability:
\begin{equation}
\ifx\du\undefined
  \newlength{\du}
\fi
\setlength{\du}{10\unitlength}

\;.
\end{equation}
The properties of $\beta$ are analogous to the properties of $\alpha$ from Subsection \ref{subsec:Inner-hom} provided that  $\Theta$ is a module trace.
\begin{proposition}
  \label{proposition:properties-beta}
  Let $\CM$ be a module category with module trace. Then the map $\beta: \Hom_{}(n,c \act m) \rr   \Hom_{}(\icm{m,n},c)$ 
    is compatible with the module structure: For all morphisms $\gamma : x \rr y$ in $\Cat{C}$ and 
    all $f \in \Hom_{}( n, c \act m)$, 
    \begin{equation}
       \beta ( \gamma \act f)= \gamma \otimes \beta(f).
    \end{equation}
\end{proposition}
\begin{proof}
   By Proposition \ref{proposition:def-beta}, $\beta(\gamma \act f)$ is uniquely determined by the requirement 
  \begin{equation}
\ifx\du\undefined
  \newlength{\du}
\fi
\setlength{\du}{10\unitlength}

 \;,
  \end{equation}
  for   all $g \in \Hom_{}(y \act (c \act m), x \act n)$. 
  From the $\Cat{C}$-compatibility of $\Theta$ and  equation (\ref{eq:specify-beta}) one obtains that  the second expression is given by
  \begin{equation}
\ifx\du\undefined
  \newlength{\du}
\fi
\setlength{\du}{10\unitlength}

\;.
  \end{equation}
  The uniqueness result of Proposition \ref{proposition:def-beta} implies that $\beta ( \gamma \act f)= \gamma \otimes \beta(f)$.
\end{proof}

Consider a module category $\CM$ with module trace.
We construct a coalgebra structure for $\icm{m,m}$ for $m \in \Cat{M}$ in analogy  to the construction of the  algebra structure of $\icm{m,m}$ in Subsection \ref{subsec:Inner-hom}.
First we define the internal coevaluation $\coev{n,m}: m \rr \icm{n,m}\act n$ as 
\begin{equation}
  \coev{n,m}= \beta^{-1} (\id_{\icm{n,m}}) \quad \corres \quad 
\ifx\du\undefined
  \newlength{\du}
\fi
\setlength{\du}{10\unitlength}

\;.
\end{equation}
Hence $\coev{n,m}$ is characterised uniquely by the property that for all $f \in \Hom_{}(\icm{n,m} \act n,m)$:
\begin{equation}
  \label{equation:def-coev}
\ifx\du\undefined
  \newlength{\du}
\fi
\setlength{\du}{10\unitlength}

\;.
\end{equation}
Next we define the internal comultiplication $\Delta_{m,n,k}:\icm{m,k} \rr \icm{n,k} \otimes 
\icm{m,n}$ by 
\begin{equation}
\ifx\du\undefined
  \newlength{\du}
\fi
\setlength{\du}{10\unitlength}

 \right)\;,
\end{equation}
and the internal counit $\epsilon: \icm{m,m} \rr 1$ as 
\begin{equation}
  \epsilon= \beta \bigg(\quad 
\ifx\du\undefined
  \newlength{\du}
\fi
\setlength{\du}{10\unitlength}
%
\;.
\end{equation}
\begin{lemma}
  \label{lemma:technical-inner-hom-coev}
  For all morphisms   $f \in \Hom_{} (\icm{n,k}\act n,c \act  m )$,
  \begin{equation}
\ifx\du\undefined
  \newlength{\du}
\fi
\setlength{\du}{10\unitlength}

\;.
  \end{equation}
\end{lemma}
\begin{proof}
  The proof is analogous to the proof of Lemma \ref{lemma:technical-inner-hom}.
\end{proof}

\begin{proposition}
  \label{proposition:inner-hom-coalgebra}
  Let $\CM$ be a module category with module trace. For any object $m \in \Cat{M}$, the internal hom $\icm{m,m}$ is canonically a coalgebra object. 
\end{proposition}
\begin{proof}
  The proof is analogous to the proof of Proposition \ref{proposition:inner-hom-algebra-obj}.
\end{proof}
It remains to prove  one more compatibility condition of $\alpha$ and $\beta$ before we can show that $\icm{m,m}$ is a Frobenius algebra.
\begin{lemma}
  \label{lemma:comp-alpha-beta}
  Consider the morphism $\coev{n,k}  \circ \ev{l,k}: \icm{l,k} \act l \rr \icm{n,k} 
  \act n$. By applying $\alpha$ and $\beta$ to this morphism we obtain the internal comultiplication and internal multiplication, respectively. In graphical terms:
  \begin{equation}
    \alpha \left( 
\ifx\du\undefined
  \newlength{\du}
\fi
\setlength{\du}{10\unitlength}

\;.
  \end{equation}
\end{lemma}
\begin{proof}
  Define $\Psi=\coev{n,k}  \circ \ev{l,k}$. First we compute $\beta(\Psi)$ using equation (\ref{eq:specify-beta}). For all $f \in \Hom_{}(\icm{n,k}) \act n, \icm{l, k} \act l)$,
  \begin{equation}
    \label{eq:beta-Psi}
\ifx\du\undefined
  \newlength{\du}
\fi
\setlength{\du}{10\unitlength}

\;, 
  \end{equation}
  where the last step involved equation (\ref{equation:def-coev}) and Lemma \ref{lemma:technical-inner-hom}. This proves that $\beta(\Psi)$ is equal to the internal multiplication. 

  Next we calculate for all $g \in \Hom_{}(\icm{n,k} \act n, \icm{l, k} \act l)$,
  \begin{equation}
\ifx\du\undefined
  \newlength{\du}
\fi
\setlength{\du}{10\unitlength}

\;, 
  \end{equation}
  where in the last step we used the definition of $\ev{l,k}$ and Lemma \ref{lemma:technical-inner-hom-coev}. Since the trace on $\Cat{C}$ is non-degenerate and $\beta$ an isomorphism, we 
  conclude that $\alpha(\Psi)$ is equal to the internal multiplication. 
\end{proof}

\begin{theorem}
  \label{theorem:Morita-to-frobenius}
  Let $\Cat{C}$ be a pivotal category and let $\Cat{M}$ be a $\Cat{C}$-module category with module trace. For all non-zero $m \in \Cat{M}$, 
  $\icm{m,m}$ is a Frobenius algebra  in $\Cat{C}$. If  $ m$ is a simple object then $\icm{m,m}$ is a
 special  haploid symmetric Frobenius algebra  with 
$\dim^{\Cat{C}}(\icm{m,m})>0$.
\end{theorem}
\begin{proof}
  We show that the relations from Definition \ref{definition:Frobenius-algebra} 
  are satisfied. Define the following morphisms for $k,l,n,r \in \Cat{M}$:
  \begin{equation}
    f_{1}= 
\ifx\du\undefined
  \newlength{\du}
\fi
\setlength{\du}{10\unitlength}

\;.
  \end{equation}
  Lemma \ref{lemma:comp-alpha-beta}, the compatibility of $\beta$ and the module action according to Proposition \ref{proposition:properties-beta} and the associativity of the  internal multiplication together imply  
  \begin{equation}
    \beta(f_{1})=\beta(f_{3}), \; \text{hence} \;f_{1}=f_{3}.
  \end{equation}
  Similarly, as a consequence of  Lemma \ref{lemma:alpha-comp-tensor}, the coassociativity of the internal comultiplication and Lemma \ref{lemma:comp-alpha-beta}, we obtain
  \begin{equation}
    \alpha(f_{2})=\alpha(f_{4}), \; \text{hence} \; f_{2}=f_{4}.
  \end{equation}
  It follows that $\alpha(f_{1})=\alpha(f_{3})$, or in graphical terms
  \begin{equation}
\ifx\du\undefined
  \newlength{\du}
\fi
\setlength{\du}{10\unitlength}

\;,
  \end{equation}
  where we again used compatibility of $\alpha$ and the module structure as well as Lemma \ref{lemma:comp-alpha-beta}.
  Similarly we conclude that  $\beta(f_{2})=\beta(f_{4})$. Together with Lemma \ref{lemma:comp-alpha-beta} and Proposition \ref{proposition:properties-beta} this implies
  \begin{equation}
\ifx\du\undefined
  \newlength{\du}
\fi
\setlength{\du}{10\unitlength}

\;.
  \end{equation}
   If we restrict attention to the case where all objects are equal to $m$, we see that  $\icm{m,m}$ satisfies   the relations (\ref{eq:Frob-equ-graphical}) defining a 
  Frobenius algebra.
  Let now $m \in \Cat{M}$ be simple. Then the identity $\Hom_{}(1,\icm{m,m}) \simeq \Hom_{}(m,m)\simeq \C$ implies that $\icm{m,m}$ is haploid. Recall that  $\eta_{m}$ and $\epsilon_{m}$ denote the internal unit and counit, 
  respectively.
  Equation (\ref{eq:specify-beta}) shows that $\epsilon_{m} \circ \eta_{m}= \Theta_{m}(\id_{m}) \neq 0$. Also by the symmetry of $\Theta$ and by equation (\ref{equation:def-coev}), 
  \begin{equation}
    \label{equation:special-frob}
\ifx\du\undefined
  \newlength{\du}
\fi
\setlength{\du}{10\unitlength}

=\dim^{\Cat{C}}(\icm{n,m}).
  \end{equation}
  As $m$ is simple,  this implies
  \begin{equation}
    \ev{n,m} \circ  \coev{n,m} = \frac{\dim^{\Cat{C}}(\icm{n,m})}{\dim^{\Theta}(m)} \cdot \id_{m}.
  \end{equation}
  Furthermore, combining Lemma \ref{lemma:comp-alpha-beta} and Lemma \ref{lemma:technical-inner-hom}, we obtain
  \begin{equation}
    \alpha \left( 
\ifx\du\undefined
  \newlength{\du}
\fi
\setlength{\du}{10\unitlength}

\;.
  \end{equation}
  Together with equation (\ref{equation:special-frob}) this implies 
  \begin{equation}
    \mu_{m,n,m}\circ \Delta_{m,n,m}=  \frac{\dim^{\Cat{C}}(\icm{n,m})}{\dim^{\Theta}(m)}  \cdot \id_{\icm{m,m}}.
  \end{equation}

  By setting $m=n$ we find that $\icm{m,m}$ is a special haploid Frobenius algebra, since by  Proposition \ref{proposition:Q-rank-1}, $\dim^{\Cat{C}}(\icm{m,m}) >0$.  
  Due to Lemma \ref{lemma:special-haploid-then-symmetric},  $\icm{m,m}$ is also a symmetric algebra.
\end{proof}

We will now prove the converse of Theorem \ref{theorem:Morita-to-frobenius}. For this  we require the following result.  An analogous statement has been proven in \cite[Lemma 2.6]{DualDef} in a slightly different setting.
\begin{lemma}
  \label{lemma:dim-modules-nonzero}
  Let $A$ be a   normalised special haploid  Frobenius algebra in a pivotal fusion category $\Cat{C}$. 
  Then  $\dim^{\Cat{C}}(M) \neq 0$ for  all simple modules $ M \in \Mod(A)$.
\end{lemma}
\begin{proof}
  The proof is a modification of the proof that all dimensions of simple objects in a pivotal fusion category are non-zero, see \cite[Lemma 2.4.1]{BakKir}.
We use the pivotal structure to identify left and right dual objects.
  First note that by Lemma \ref{lemma:special-haploid-then-symmetric}, $A$ is symmetric and  for a symmetric Frobenius algebra, 
  \begin{equation}  
 \label{symmetric-then-inverse}
\ifx\du\undefined
  \newlength{\du}
\fi
\setlength{\du}{10\unitlength}
 \;.
  \end{equation}
  This follows from the fact that the left hand side is the inverse of the  morphism on the left of equation (\ref{eq:symmetric-morphisms}), while the right hand side is the inverse morphism of 
  the right hand side of (\ref{eq:symmetric-morphisms}), hence both have to agree.

  Let $(M, \rho)$ be a simple $A$-module. 
  Proposition \ref{proposition:modules-non-zero} implies 
  $\C=\Hom_{A}(M,M) \simeq \Hom_{\Cat{C}}( M \tensor[A] {}^{*}M,1)$. 
  It is sufficient to show that there are non-zero maps 
  $\coev{M}^{A}: 1 \rr M \tensor[A] {}^{*}M$ and $\ev{M}^{A}:M \tensor[A] {}^{*}M \rr 1$ for which the diagram
  \begin{equation}
    \label{eq:dim-simple-M-neq-0}
    \begin{tikzcd}
              & 1 \ar{dl}[above, xshift=-22]{(1 \otimes {}^{*}a_{M})\coev{M}} \ar{d}{\coev{M}^{A}} \\
        M \otimes {}^{*}M \ar{r}{P} \ar{dr}[below]{\evp{M}} &M \tensor[A]{}^{*}M \ar{d}{\ev{M}^{A}}\\
        & 1
         \end{tikzcd}
  \end{equation}
  commutes. The semisimplicity of $\Cat{C}$ then implies that the composition $\ev{M} \circ \coev{M}$ is non-zero.
  We obtain $\ev{M}^{A}$ from the universal property of the cokernel by observing that 
$\evp{M} \circ ( \rho \otimes \id_{{}^{*}M})= \ev{M} \circ ( \id_{M} \otimes \rho_{{}^{*}M})$ as morphisms 
  $M \otimes A \otimes {}^{*}M \rr 1$.  Here
  $\rho_{{}^{*}M}$ is defined by (\ref{equation:dual-action}). For $\coev{M}^{A}$ we have to show that 
$P \circ (1 \otimes {}^{*}a_{M}) \circ  \coev{M} \neq 0$, where $P$ is the projector (\ref{equation:projector1}). 
We compute, using in the third step the pivotal structure to turn the left dual of the action $\rho_{M}$ of $A$ on $M$
into the right dual $\rho_{M}^{*}$:
   \begin{align} 
    &P \circ (1 \otimes {}^{*}a_{M}) \circ \coev{M}  \corres 
\ifx\du\undefined
  \newlength{\du}
\fi
\setlength{\du}{10\unitlength}
=(1 \otimes {}^{*}a_{M}) \circ \coev{M}.
  \end{align}
 This proves  the statement.  
\end{proof}
\begin{proposition}
  \label{proposition:converse-statement-frob-mod}
  Let $A$ be a special haploid symmetric  Frobenius algebra in $\Cat{C}$. Then the $\Cat{C}$-module category of 
  right $A$-modules, $\Mod(A)$, has a module trace induced by the trace on $\Cat{C}$. In particular, $A$ satisfies $\dim^{\Cat{C}}(A) >0$.
\end{proposition}
\begin{proof}
  The symmetry and $\Cat{C}$-compatibility follow from the properties of the trace $\tr^{\Cat{C}}$ in $\Cat{C}$. We only have to show that the induced pairing on the $\Hom_{}$-spaces of  $\Mod(A)$ is non-degenerate. 
  According to Lemma \ref{lemma:general-trace-ses} it is sufficient to show that all simple modules $m$ over $A$ have $\dim^{\Cat{C}}(m) \neq 0$.
  This follows from Lemma \ref{lemma:dim-modules-nonzero}. 

Consider the quantum dimension of $A$. Since $A$ is haploid it is a simple module over itself. The inner hom object of $\Mod(A)$ is given by 
the tensor product over $A$, hence $\icm{A,A} = A \tensor[A] A=A$, see e.g. \cite{SchwCat} for the last equality.  The statement  now follows from Proposition \ref{proposition:Q-rank-1}.
\end{proof}
We have thus established the correspondence between module traces and Frobenius algebras. 

The following example illustrates  the role of the pivotal structure of $\Cat{C}$ in this correspondence. 
Consider $\Cat{C}= \Vect_{G}$ be the category of $G= \Z / 2 \Z$-graded vector spaces. 
The sum of the two simple objects in $\Cat{C}$ defines an algebra $A$. If we choose the non-standard pivotal structure for $\Cat{C}$, where the simple object corresponding to $-1 \in G$ has dimension $-1$, then $A$ has dimension zero and cannot possess the structure of a  special symmetric Frobenius algebra. This agrees with the discussion 
in Example \ref{example:finite-G}, which implies that the corresponding module category does not possess a module trace.

If $\CM$ is a module category with module trace, the dimensions of simple objects in general change under the equivalence 
$ \CM \ni n \mapsto \icm{m,n} \in \Mod(\icm{m,m})$  with $m \in \Cat{M}$ a simple object. The following lemma allows one to calculate the relevant scaling factor.
\begin{lemma}
  Let $\CM$ be a module category with module trace. Let $m,n \in \Cat{M}$ be simple objects. Then 
  \begin{equation}
    \label{eq:change-of-dimension}
    \dim^{\Cat{C}}(\icm{m,n}) =\frac{\dim^{\Cat{C}}(\icm{m,m})}{\dim^{\Theta}(m)} \cdot \dim^{\Theta}(n).
  \end{equation}
  Under the equivalence $\CM \simeq \Mod(\icm{m,m})$ the dimensions of simple objects 
  are scaled by $\frac{\dim^{\Cat{C}}(\icm{m,m})}{\dim^{\Theta}(m)}$.
\end{lemma}

\begin{proof}
  Set $d_{i}= \dim^{\Theta}(m_{i})$. From Proposition \ref{proposition:Q-rank-1}, we obtain
  \begin{equation}
    \label{eq:consequence-rank-1}
    Q_{ij}= \frac{d_{i}\cc{d_{j}} \dim(\Cat{C})}{\sum_{k}|d_{k}|^{2}}. 
  \end{equation}
  This implies
  \begin{equation}
    \begin{split}
      \dim^{\Cat{C}}(\icm{m_{j},m_{i}})&=d_{i} \cdot \frac{|d_{j}|^{2} \dim(\Cat{C})}{d_{j} \sum_{k}|d_{k}|^{2}}\\
      &=d_{i} \cdot \frac{\dim^{\Cat{C}}(\icm{m_{j},m_{j}})}{d_{j}},
    \end{split}
  \end{equation}
 where we again used equation (\ref{eq:consequence-rank-1}) in the last step. Setting
  $m=m_{j}$ and $n=m_{i}$ then proves the claim.
\end{proof}
Finally we interpret our result using the notion of Morita equivalence of algebras  (see \cite{Ostrik}).
Two algebras $A,B \in \Cat{C}$ are called Morita equivalent if the categories $\Mod(A)$ and $\Mod(B)$ are equivalent as module categories. 
\begin{theorem}
  \label{theorem:separable-Morita}
  Every separable indecomposable algebra $A$ in a  fusion category with a flexible pivotal structure is Morita equivalent to a special haploid symmetric Frobenius algebra. 
\end{theorem}
\begin{proof}
  By definition of a flexible pivotal structure, the module category $\Mod(A)$ possesses a module trace. By Theorem \ref{theorem:Morita-to-frobenius}, this module category is equivalent to the module category
  corresponding to a special haploid Frobenius algebra. 
\end{proof}
Together with Proposition \ref{proposition:pseudo-unitary} this implies the following:
\begin{corollary}
 If an indecomposable algebra $A$ in a pseudo-unitary fusion category  $\Cat{C}$ is separable, then it is Morita equivalent to a  special haploid symmetric Frobenius algebra. 
\end{corollary}

\appendix
\section{Graphical Calculus for Tensor Categories}
\label{subsection:Graphical-calculus}
We summarise the  graphical calculus for tensor categories, see e.g. \cite{BakKir}.
The symbol $\corres$ is used to indicate that  a certain diagrammatic expression represents 
an algebraic expression. Objects in $\Cat{C}$ and the tensor product are represented by the following diagrams.
\begin{equation}
  c \corres 
\ifx\du\undefined
  \newlength{\du}
\fi
\setlength{\du}{10\unitlength}

.
\end{equation}
Morphisms are  represented by labelled boxes, and we do not distinguish objects from their unit morphisms.
All diagrams are read from top to bottom. The  composition is given by vertical connection of boxes.
\begin{equation}
  f: c \rr d \quad \corres \quad 
\ifx\du\undefined
  \newlength{\du}
\fi
\setlength{\du}{10\unitlength}

.
\end{equation}
An empty box represents a $\Hom_{}$-vector space:  
\begin{equation}
  \Hom_{}(c,d) \quad \corres \quad 
\ifx\du\undefined
  \newlength{\du}
\fi
\setlength{\du}{10\unitlength}
%
\;.
\end{equation}
The tensor product of two morphisms $f: c \rr d$ and $g: a \rr b$ is depicted as follows:
\begin{equation}
  f \otimes g \quad \corres \quad 
\ifx\du\undefined
  \newlength{\du}
\fi
\setlength{\du}{10\unitlength}

.
\end{equation}
The interchange law $f \otimes g= (f \otimes \id_{a})(\id_{d} \otimes g)= (\id_{c} \otimes g)(f \otimes \id_{b})$ has the following graphical expression:
\begin{equation}
\ifx\du\undefined
  \newlength{\du}
\fi
\setlength{\du}{10\unitlength}

.
\end{equation}
The graphical notation suppresses the unit object and the associativity constraint in $\Cat{C}$. Due to Mac Lane's coherence theorem for monoidal categories, 
a graphical expression uniquely defines a  morphisms in $\Cat{C}$ once parentheses 
and unit objects are specified for the incoming and outgoing objects. 
The evaluation and coevaluation morphisms for the right duals are depicted as follows: 
\begin{equation}
  \ev{c} \quad \corres \quad 
\ifx\du\undefined
  \newlength{\du}
\fi
\setlength{\du}{10\unitlength}

,
\end{equation}
and the rigidity axioms read:
\begin{equation}
  \label{rigidity}
\ifx\du\undefined
  \newlength{\du}
\fi
\setlength{\du}{10\unitlength}
%
.
\end{equation}
The graphical notation for left duals is analogous. 
If $\Cat{C}$ is a pivotal category, the pivotal isomorphism is represented graphically by a smaller box, 
which is often not labelled.  

The right dual of a morphism $f:c \rr d$ is defined by: 
\begin{equation}
\ifx\du\undefined
  \newlength{\du}
\fi
\setlength{\du}{10\unitlength}

.
\end{equation}
The left dual of a morphism is defined analogously using the left duality, and the map 
\begin{equation}
\ifx\du\undefined
  \newlength{\du}
\fi
\setlength{\du}{10\unitlength}
%
\end{equation}
yields an isomorphism $\Hom_{}(c \otimes d, x) \rr \Hom_{} (d, \leftidx{^{*}}{c}{} \otimes x)$.

\subsection*{Acknowledgements}
The author thanks Catherine Meusburger for  many helpful discussions, support and detailed reading of the draft. He furthermore thanks 
Christoph Schweigert,  Martin Mombelli and Ingo Runkel for their hospitality, useful discussions and support, Martina Lanini, Winston Fairbairn and Nils Carqueville for proofreading and 
Torsten Sch\"onfeld  for his help with Latex and the implementation of the diagrams. 
We thank Pavel Etingof for sharing his answers to the questions in Section \ref{sec:module-traces}
and for interesting discussions. 
 The author is grateful to the referee for useful comments.

This work was funded by the research grant ME 3425/1-1 in the Emmy-Noether program of the German Research Foundation.

\bibliographystyle{hplain}	
\bibliography{bimod}

\end{document}